\documentclass{icmart}
\usepackage{graphicx,url}
\usepackage[numbers,sort]{natbib}
\newcommand{\hf}{\mbox{$\tfrac{1}{2}$}}

\contact[imj@stanford.edu]{Department of Statistics, Sequoia Hall, Stanford University, Stanford CA 94305, U.S.A.}

\theoremstyle{definition}

\title[High Dimensional Statistical Inference and Random Matrices]{High Dimensional Statistical Inference and Random Matrices}

\author[Iain M. Johnstone]{Iain M. Johnstone\thanks{The author is
    grateful to Persi Diaconis, Noureddine El Karoui, Peter Forrester,
    Matthew Harding, Plamen Koev, Debashis Paul, Donald Richards and
    Craig Tracy for advice and comments during the writing of this
    paper, to the Australian National University for hospitality, and
    to NSF DMS 0505303 and  NIH R01 EB001988 for financial support.}}

\begin{document}

\begin{abstract}
Multivariate statistical analysis is concerned with observations on
several variables which are thought to possess some degree of
inter-dependence. Driven by problems in genetics and the social
sciences, it first flowered in the earlier half of the last century.
Subsequently, random matrix theory (RMT) developed, initially within
physics, and more recently widely in mathematics.  While some of the
central objects of study in RMT are identical to those of multivariate
statistics, statistical theory was slow to exploit the connection.
However, with vast data collection ever more common, data sets
now often have as many or more variables than the number of individuals
observed.  In such contexts, the techniques and results of RMT have
much to offer multivariate statistics.
The paper reviews some
of the progress to date.
\end{abstract}

\begin{classification}
Primary 62H10; 62H25; 62H20; Secondary 15A52.
\end{classification}

\begin{keywords}
canonical correlations; 
eigenvector estimation;
largest eigenvalue;
principal components analysis;
Random matrix theory; 
Wishart distribution;
Tracy-Widom distribution.
\end{keywords}

\maketitle

\section{Introduction}

Much current research in statistics, both in statistical theory, and
in many areas of application, such as genomics, climatology or
astronomy, focuses on the problems and opportunities posed by
availability of large amounts of data. (More detail may be found, for
example, in the paper by Fan and Li \cite{fan06} in these proceedings.)
There might be many variables and/or many observations on each
variable.  Loosely one can think of each variable as an additional
dimension, and so many variables corresponds to data sitting in a high
dimensional space.
Among several mathematical themes one could follow -- Banach space
theory, convex geometry, even topology -- this paper focuses on Random
Matrix Theory, and some of its interactions with important areas of what
in statistics is called ``Multivariate Analysis.''

Multivariate analysis deals with observations on more
than one variable when there is or may be some dependence between the
variables.  The most basic phenomenon is that of correlation -- the
tendency of quantities to vary together: tall parents tend to have
tall children. 
From the beginning, there has also been a focus on summarizing and
interpreting data by reducing dimension, for example by methods such
as Principal Components Analysis (PCA).
While there are many methods and corresponding problems of
mathematical interest, this paper concentrates largely
on PCA as a leading example, together with a few remarks on related
problems. 
Other overviews with substantial statistical content include 
\cite{bai99}, \cite{diac03} and~\cite{edra05}.

In an effort to define terms and give an example, the earlier sections
cover introductory material, to set the stage.  The more recent work,
in the later sections, concentrates on results and phenomena which
appear in an asymptotic regime in which $p$, the number of variables
increases to infinity, in proportion to sample size~$n$.

\section{Background}
\label{sec:Background}

\subsection{Principal Components Analysis}
\label{sec:PCA}

Principal Components Analysis (PCA) is a standard technique of
multivariate statistics, going back to Karl Pearson in 1901
\cite{pear01} and Harold Hotelling in 1933 \cite{hote33}. There is a
huge literature \cite{joll02} and interesting modern variants continue
to appear \cite{tsl00,rosa00}.
 A brief description of the classical method,
an example and references are included here for convenience.

PCA is usually described first for abstract random variables, and then
later as an algorithm for observed data.  So first suppose we have $p$
variables $\mathsf{X}_1, \ldots, \mathsf{X}_p$.  We think of these as
random variables though, initially, little more is assumed than the
existence of a \textit{covariance matrix} $\Sigma = (\sigma_{k k'})$,
composed of the mean-corrected second moments
\begin{displaymath}
  \sigma_{k k'} = \text{Cov}(\mathsf{X}_k, \mathsf{X}_{k'}) =
  E(\mathsf{X}_k - \mu_k)(\mathsf{X}_{k'}-\mu_{k'}).
\end{displaymath}

The goal is to reduce dimensionality by constructing a smaller number
of ``derived'' variables $W = \sum_k v_k \mathsf{X}_k$, having variance
\begin{displaymath}
 \text{Var}(W) = \sum_{k,k'} v_k \sigma_{kk'} v_{k'} = \mathbf{v}^T
 \Sigma \mathbf{v}.   
\end{displaymath}
To concentrate the variation in as few derived variables as possible,
one looks for vectors that maximize $\text{Var}(W)$.
Successive linear combinations are sought that are
orthogonal to those previously chosen.  The \textit{principal
  component eigenvalues} $\ell_j$ and \textit{principal component
  eigenvectors} $\mathbf{v}_j$ are thus obtained from
\begin{equation}
\label{eq:popn}
  \ell_j = \max \{ \mathbf{v}^T \Sigma \mathbf{v} ~:~ \mathbf{v}^T
\mathbf{v}_{j'} = 0 ;~ j' < j, \ |\mathbf{v}|=1  \}.
\end{equation}

In statistics, it is common to assume a stochastic model in terms of
random variables whose distributions contain unknown parameters, which
in the present case would be the covariance matrix and its resulting
principal components.  To \textit{estimate} the unknown parameters of
this model we have observed data, assumed to be $n$ observations on
each of the $p$ variables. The observed data on variable
$\mathsf{X}_k$ is viewed as a vector $\mathbf{x}_k \in \mathbb{R}^n$.
The vectors of observations on each variable are collected as rows
into a $p \times n $ data matrix
\begin{displaymath}
  X = (x_{ki}) = [\mathbf{x}_1 \ldots \mathbf{x}_p]^T.
\end{displaymath}

\begin{figure}
  \centering
  \centerline{\includegraphics[width=0.85\textwidth]{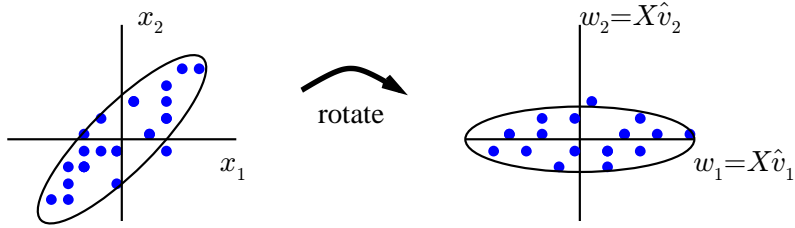}}

  \caption{The $n$ data observations are viewed as $n$ points in $p$
    dimensional space, the $p$ dimensions corresponding to the
    variables. The sample PC eigenvectors $\hat{\mathbf{v}}_j$ create
    a rotation of the variables into the new derived variables, with
    most of the variation on the low dimension numbers. In this two
    dimensional picture, we might keep the first dimension and discard
    the second.}
  \label{fig:pceg}
\end{figure}

A standard pre-processing step is to center each variable by subtracting
the sample mean $\bar x_k = n^{-1} \sum_i x_{ki}$, so that $x_{ki}
\leftarrow x_{ki} - \bar x_k.$ After this centering, define the
$p \times p$ \textit{sample covariance matrix} $S = (s_{kk'})$ by
\begin{displaymath}
 S = (s_{k k'}) = n^{-1} X X^T, \qquad 
   s_{k k'} = n^{-1} \sum_i x_{ki}x_{k'i}.
\end{displaymath}
The derived variables in the sample, $\mathbf{w} = X \mathbf{v} = \sum_k v_k
\mathbf{x}_k,$ have sample variance $ \widehat{\text{Var}}(\mathbf{w})
= \mathbf{v}^T S \mathbf{v}.$    
Maximising this quadratic form leads to successive sample principal
components $\hat{\ell_j}$ and $\hat{\mathbf{v}}_j$ from
the sample analog of (\ref{eq:popn}):
\begin{displaymath}
     \hat \ell_j = \max \{ \mathbf{v}^T S \mathbf{v} ~:~ \mathbf{v}^T
      \hat{\mathbf{v}}_{j'} = 0, \ j' < j, \ |\mathbf{v}| = 1 \} 
\end{displaymath}
Equivalently, we obtain for $j = 1, \ldots, p$,
\begin{displaymath}
   S \hat{\mathbf{v}}_j = \hat \ell_j \hat{\mathbf{v}}_j, \qquad
   \hat w_j = X \hat{\mathbf{v}}_j.
\end{displaymath}
Note the statistical convention: estimators derived from samples
are shown with hats. 
Figure \ref{fig:pceg} shows a conventional picture illustrating PCA.

Observed data are typically noisy, variable, and limited in quantity,
so we are interested in the estimation errors
\begin{displaymath}
 \hat \ell_j(X) -  \ell_j, \qquad \qquad  \widehat{\mathbf{v}}_j(X)  - 
 \mathbf{v}_j.
\end{displaymath}
An additional key question in practice is: how many dimensions are
``significant'', or should be retained? One standard approach is to
look at the percent of total variance explained by each of the principal
components:
\begin{displaymath}
  p_j = \hat \ell_j/ \sum \hat \ell_{j'} = 
      \hat \ell_j/ \text{tr} S.
\end{displaymath}

\textit{An example.} \  \citet*{mpcs78}
is a celebrated example of the use of PCA in human
genetics and anthropology.  It was known from archaeological
excavations that farming spread gradually from Near East across Europe
9000--5000 yrs ago (map below right).  A motivating question was whether
this represented spreading of the farmers themselves (and hence their
genes) or transfer of technology to pre-existing populations (without
a transfer of genes).

\citet{mpcs78} brought genetic data to bear on the issue.
Simplifying considerably, the data matrix $X$ consisted of
observations on the frequencies of alleles of $p = 38$ genes in human
populations at $n = 400$ locations in Europe.  The authors sought to
combine information from the 38 genes to arrive at a low dimensional
summary.

A special feature of the genetics data is that the observations $i$
have associated locations $\text{loc}[i]$, so that it is possible to
create a map from each of the principal components $w_j$, by making a
contour plot of the values of the derived variable $w_j[i]$ at each of
the sampling locations $\text{loc}[i]$.  For example the first
principal component (map below left) shows a clear trend from south-east
to north-west, from Asia Minor to Britain and Scandinavia.  The
remarkable similarity of the PC map, derived from the gene
frequencies, with the farming map, derived from archaeology, has been
taken as strong support for the spread of the farmers themselves.

For the genetics data, the first component (out of 38) explains $p_1 =
27\%$ of the variance, the second $p_2 = 18\%$, and the third $p_3 =
11\%$. Thus, and this is typical, more than half the variation is
captured in the first three PCs.  The second and third, and even
subsequent PCs also show patterns with important linguistic and
migratory interpretations. For more detail, we refer to books of
Cavalli-Sforza \citep{csmp94,cs00}, from which the maps below are reproduced.

\bigskip

\parbox{.55\textwidth}{
\centerline{\includegraphics[width = .4\textwidth, angle = 90]{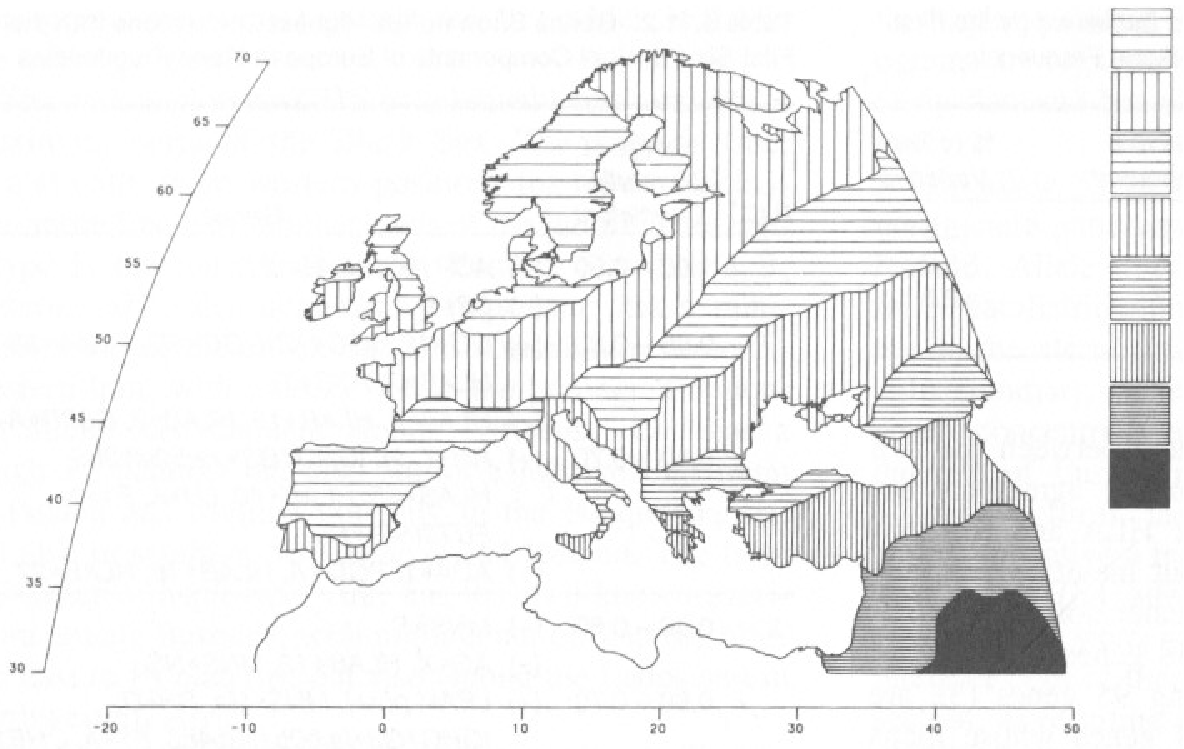}}
}
\parbox{.44\textwidth}{
\centerline{\includegraphics[width = .45\textwidth, angle = 0]{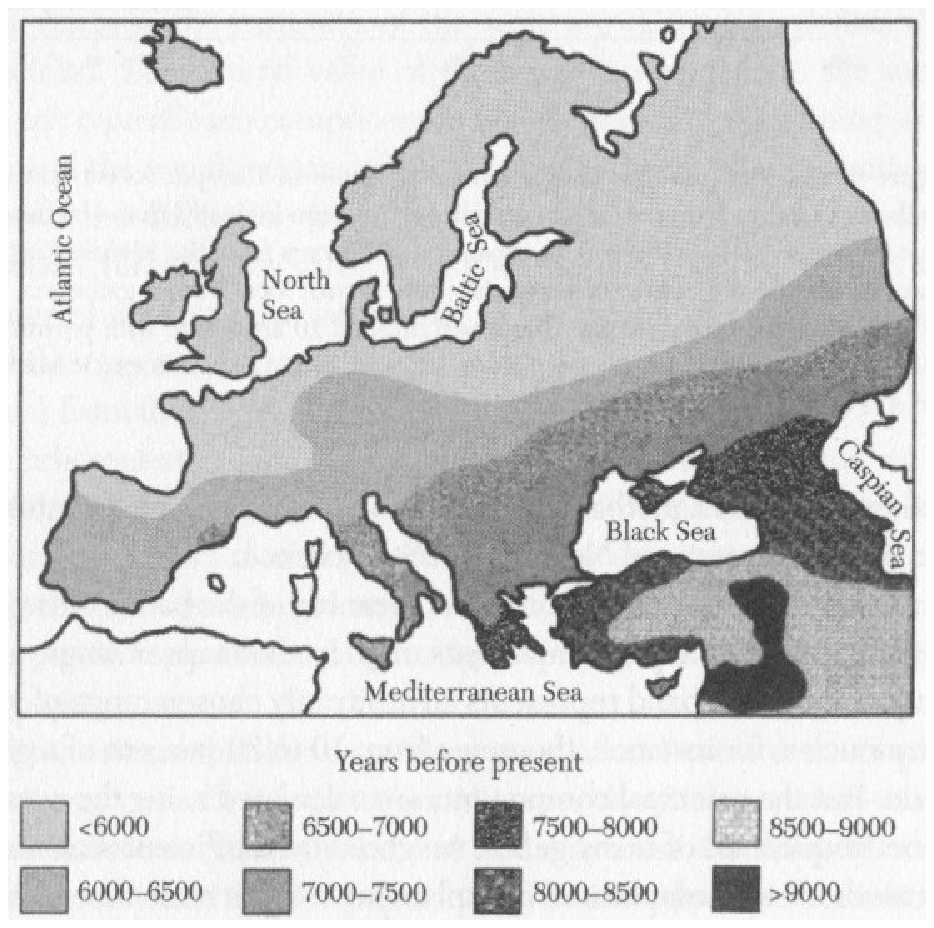}}
}

 \bigskip

\subsection{Gaussian \& Wishart Distributions}

For quantitative analysis, we need  more specific assumptions
about the process generating the data.  The simplest and most
conventional model assumes that the $p$ random variables
$\mathsf{X}_1, \ldots, \mathsf{X}_p$ follow a $p-$variate Gaussian distribution
$N_p(\mu, \Sigma)$, with mean $\mu$ and covariance matrix $\Sigma$,
and with probability density function for $\mathsf{X} =
(\mathsf{X}_1, \ldots, \mathsf{X}_p)$ given by
\begin{displaymath}
f(\mathsf{X}) = |\sqrt{2 \pi}
\Sigma|^{-1/2} \exp \{ -\tfrac{1}{2} (\mathsf{X}-\mu)^T \Sigma^{-1}
(\mathsf{X}-\mu) 
\}.  
\end{displaymath}

The observed sample is assumed to consist of $n$ independent draws
$X_1, \ldots, X_n$ from $\mathsf{X} \sim N_p(\mu, \Sigma),$
collected into a $p \times n$ data matrix $X = [X_1 \ldots X_n]$. 
When focusing on covariances, it is a slight simplification to assume
that $\mu = 0$, as we shall here. 
In practice, this idealized model of
independent draws from a Gaussian is generally at best approximately
true -- but we may find some reassurance in the dictum ``All models are
wrong, some are useful.'' \cite{box79}

The (un-normalized) cross product matrix 
$A = X X^T$ is said to have a $p$ - variate \textit{Wishart}
distribution on $n$ degrees of freedom. The distribution is named for
John Wishart who in 1928 \cite{wish28} derived the density function
\begin{displaymath}
  f(A) = c_{n,p} | \Sigma |^{-n/2} |A|^{(n-p-1)/2} \exp \{ -
  \tfrac{1}{2} \text{tr} (\Sigma^{-1} A) \},
\end{displaymath}
which is supported on the cone of non-negative definite matrices.
Here $c_{n,p}$ is a normalizing constant, and it is assumed that $\Sigma$
is positive definite and that $n \geq p$.  

The eigendecomposition of the Wishart matrix connects directly with
Principal Components Analysis.
Start with a Gaussian data matrix, form the covariance $S$, yielding a
Wishart density for $A = nS$.  The eigenvalues and vectors of $A$,
given by
\begin{equation}
  \label{eq:single-w}
  A u_i = l_i u_i,  \qquad l_1 \geq \ldots \geq l_p \geq 0,
\end{equation}
are related to the principal component eigenvalues and vectors \textit{via}
\begin{displaymath}
    l_i = n \hat \ell_i, \qquad  u_i = \hat{\mathbf{v}}_i. 
\end{displaymath}

\medskip
\textit{Canonical Correlations.} \ 
We digress briefly from the PCA theme to mention one additional multivariate
technique, also due to Hotelling \citep{hote36}, since it will help
indicate the scope of the results. 
Given two sets of variables $\mathsf{X} = ( \mathsf{X}_1, \ldots,
\mathsf{X}_p)$ and $\mathsf{Y} = ( \mathsf{Y}_1, \ldots,
\mathsf{Y}_q)$, with
a joint $p+q$-variate Gaussian distribution,
we may ask for that linear combination of $\mathsf{X}$ that is
most correlated with some linear combination of $\mathsf{Y}$, seeking
the canonical correlations
\begin{displaymath}
  r_i^2 = \max_{u_i,v_i} \mbox{Corr} \, (u^T_i \mathsf{X}, v^T_i \mathsf{Y}),
\end{displaymath}
and the maximization is subject to $|u_i| = |v_i| = 1$.
 
To take an example from climatology \cite{bapr87}: the $\mathsf{X}$
variables might be sea surface temperatures at various ocean
locations, and the $\mathsf{Y}$ variables might be land temperatures
at various North American sites. The goal may be to find the
combination of sea temperatures that is most tightly correlated with
some combination of land temperatures.
For a recent example in functional magnetic resonance imaging, see 
\cite{frim01}.

If we have $n$ draws $(X_i, Y_i), \ i=1, \ldots, n$ from
the joint distribution, the
sample version of this problem may be written as a generalized
eigenequation that involves \textit{two} independent matrices $A$ and
$B$, each following $p-$variate Wishart distributions -- on $q$ and
$n-q$ degrees of freedom respectively: 
\begin{displaymath}
 A v_j = r_j^2 (A+B) v_j, \qquad r_1^2 \geq \ldots \geq r_p^2.  
\end{displaymath}

The parameters of the Wishart distribution depend on those of the
parent Gaussian distribution of the data -- if $X$ and $Y$ are
independent, then they both reduce to Wishart matrices with identity
covariance matrix: $A \sim W_p(q,I)$ and $B \sim W_p(n-q,I)$.

\medskip

\textit{The Double Wishart setting.} \ 
Suppose we have two independent Wishart matrices $A \sim W_p(n_1,I)$
and $B \sim W_p(n_2,I)$, with the degrees of freedom parameters $n_1,
n_2 \geq p.$ We call this the double Wishart setting.  Two remarks: By
writing Wishart distributions with \textit{identity} matrices, we
emphasize, for now, the ``null hypothesis'' situation in which there
is no assumed structure (compare Section \ref{sec:beyond}).  Second,
by taking a limit with $n_2 \to \infty$, one recovers the single
Wishart setting.

Of central interest are the roots $x_i, i = 1, \dots, p$ of the generalized
eigenproblem constructed from $A$ and $B$: 
\begin{equation}
  \label{eq:double-w}
  \det [ x(A+B) - A ] = 0.
\end{equation}

The canonical correlations problem is a leading example.  In addition,
essentially all of the classical multivariate techniques involve an
eigendecomposition that reduces to some form of this equation. 
Indeed, we may collect almost all the chapter titles in any classical
multivariate statistics textbook
(e.g. \cite{ande84,muir82,mkb79,jowi02}) into a table: 

\begin{table}[h]
  \centering
\begin{tabular}{ll}
\textsf{Double Wishart}  & \textsf{Single Wishart} \\
Canonical correlation analysis &  Principal Component analysis \\
Multivariate Analysis of Variance 
&  Factor analysis \\
Multivariate regression analysis &  Multidimensional scaling  \\
Discriminant analysis & \\
Tests of equality of covariance matrices &
\end{tabular}
\end{table}

This table emphasizes the importance of finding the distribution of
the roots of (\ref{eq:double-w}), which are 
basic to the use of these methods in applications.

\bigskip

\textit{Joint density of the eigenvalues.} \ The joint null hypothesis
distribution of the eigenvalues for canonical correlations and
principal components was found in 1939. The results were more or less
simultaneously obtained by five distinguished statisticians in three
continents \citep{fish39,girs39,hsu39,mood51,roy39}:
\begin{equation} \label{eq:jtdens}
  f(x_1, \ldots, x_p) = c \prod_i w^{1/2}(x_i) \, \prod_{i<j} (x_i - x_j),
\qquad \qquad x_1 \geq \ldots \geq x_p,
\end{equation}
with
\begin{displaymath}
  w(x) =
  \begin{cases}
    x^{n-p-1} e^{-x}  & \text{single Wishart} \\
    x^{n_1 - p -1} (1-x)^{n_2-p-1} & \text{double Wishart}.
  \end{cases}
\end{displaymath}
The normalizing constant $c$ is given, using the multivariate Gamma function
  $\Gamma_p(a) = \pi^{p(p-1)/4} \prod_{i=1}^p \Gamma(a - (i-1)/2)$, by
\begin{displaymath}
 c =
 \begin{cases}
   \frac{2^{-{pn/2}}  \pi^{p^2/2}}{ \Gamma_p(p/2) \Gamma_p(n/2)} & 
 \text{single Wishart} \\
   \frac{ \pi^{p^2/2} \Gamma_p((n_1+n_2)/2)}{\Gamma_p(p/2) \Gamma_p(n_1/2)
   \Gamma_p(n_2/2)} & \text{double Wishart}.
 \end{cases}
\end{displaymath}

Thus, the density  has a product term involving each of the roots one at
a time, through a weight function $w$ which one recognizes as the
weight function for two of the classical families of orthogonal
polynomials, Laguerre and Jacobi respectively.

The second product is the so-called ``Jacobian'' term, which arises in
the change of variables to eigenvalue and eigenvector co-ordinates. It is
this pairwise interaction term, also recognizable as a Vandermonde
determinant (see (\ref{eq:vandermonde}) below), that causes difficulty in the
distribution theory.

This result was the beginning of a rich era of multivariate
distribution theory in India, Britain, the U.S., and Australia,
summarized, for example, in \cite{ande84,muir82,mkb79}. While some of
this theory became so complicated that it lost much influence on
statistical practice, with new computational tools and theoretical
perspectives the situation may change.

\subsection{Random Matrices}
\label{sec:RMT}

We detour around this theory and digress a moment to introduce
the role of random matrix theory.
Beginning in the 1950s, physicists began to use random matrix models to
study quantum phenomena. 
In quantum mechanics the energy levels of a system, such as the
nucleus of a complex atom, are described by the eigenvalues of a
Hermitian operator $H$, the Hamiltonian:
$H \psi_i = E_i \psi_i$, with $E_0 \leq E_1 \leq \cdots $.
The low-lying energy levels can be understood by theoretical work, but
at higher energy levels, for example in the millions, the
analysis becomes too complicated.

Wigner proposed taking the opposite approach, and sought a purely
statistical description of an ``ensemble" of energy levels -- that
could yield properties such as their empirical distribution and the
distribution of spacings. He further made the hypothesis that the
local statistical behavior of energy levels (or eigenvalues) is well
modeled by that of the eigenvalues of a random matrix.  Thus the
approximation is to replace the Hermitian operator $H$ by a large
finite \textit{random} $N \times N$ matrix $H_N$.

One example of a statistical description that we will return to later
is the celebrated SemiCircle Law \cite{wign55,wign58}.  This refers to the
eigenvalues of a so-called Wigner matrix $H_N$, with independent and
identically distributed entries of mean $0$ and a finite
variance $\sigma^2$. With no further conditions on the distribution of
the matrix entries, the empirical distribution 
$F_N(t) = \# \{ i : x_i \leq t \}/N$
of the eigenvalues
converges to a limit with density given by a semicircle:
\begin{displaymath}
 dF_N(x \sigma \sqrt{N}) \rightarrow \frac{1}{4\pi} \sqrt{ 4 - x^2 } dx.
\end{displaymath}

\medskip \textbf{ Ensembles and Orthogonal Polynomials } \ 
Quite early on, there was interest in eigenvalue distributions whose
densities could be described by more general families of weight
functions than the Gaussian.  For example, \citet{foka64} used weight
functions from the families of classical orthogonal polynomials.
Analogies with statistical mechanics made it natural to introduce an
additional (inverse temperature) parameter $\beta$, so that the
eigenvalue density takes the form
\begin{equation}
  \label{eq:fox-kahn}
  f(x_1, \ldots, x_N) = 
  c \prod_1^N w(x_i)^{\beta/2} \prod_{i<j} |x_i - x_j|^\beta.
\end{equation}

At this time, it was only partially realized that in the case
$\beta= 1$, these densities were already known in statistics. 
But the table shows that in fact, the three classical orthogonal polynomial
weight functions correspond to the 
three most basic null eigenvalue
distributions in multivariate statistics:

\begin{table}[h]
  \centering
\begin{tabular}{lllll}
  $w(x) =$  &  $  e^{-x^2/2}$  &  Hermite  & $H_k$ &    \textit{Gaussian} \\
            &  $ x^a e^{-x}$  &  Laguerre & $L^a_k$ &   \textit{Wishart}  \\
    &  $(1-x)^a (1+x)^b$  &  Jacobi  & $P^{a,b}_k$ &  \textit{Double Wishart} 
\end{tabular}  
  \caption{The orthogonal polynomials are taken in the standard forms
    given in \citet{szeg67}.}
  \label{tab:polys}
\end{table}

\citet{dyso62} showed that physically reasonable symmetry assumptions
restricted $\beta$ to one of three values: 
\begin{center}
\begin{tabular}{lll}
               &   \textsf{Symmetry Type}   &  \textsf{Matrix entries}  \\
  $\beta = 1$  &   orthogonal  &  real            \\
  $\beta = 2$  &   unitary     &  complex         \\
  $\beta = 4$  &   symplectic  &  quaternion 
\end{tabular}  
\end{center}
Mathematically, the complex-valued case is always
the easiest to deal with, but of course it is the real case that is of
primary statistical (and physical) interest; though cases with complex
data do occur in applications, notably in communications. 

To summarize, the classical ``null hypothesis'' distributions in
multivariate statistics correspond to the \textit{italicized}
eigenvalue densities in the

\smallskip
\centerline{$\begin{Bmatrix}
  \text{Gaussian} \\  \text{\textit{Laguerre}} \\  \text{\textit{Jacobi}}
\end{Bmatrix}
\begin{Bmatrix}
  \text{\textit{Orthogonal}} \\  \text{\textit{Unitary}} \\ \text{Symplectic}
\end{Bmatrix}$
Ensemble.}

\smallskip
\noindent These are often abbreviated to LOE, JUE, etc. We have not
italicized the Symplectic case for lack (so far) of motivating
statistical applications (though see \cite{abj83}).

\bigskip \textbf{ Some uses of RMT in Statistics } \ 
This table organizes some of the classical topics within RMT, and some
of their uses in statistics and allied fields. 
This paper will focus selectively (topics in italics), and
in particular on largest eigenvalue results and their use for an 
important class of hypothesis tests, where RMT brings something quite
new in the approximations.

\begin{table}[h]
  \centering
\begin{tabular}{ll}
  \textit{Bulk}  &  Graphical methods \cite{wach78,wach80} [finance
  \cite{bopo03,pbl05}, communications \cite{tuve04}] \\
    Linear Statistics & Hypothesis tests, distribution theory \\
    \textit{Extremes} &  \textit{Hypothesis tests, distribution theory,}
  role in proofs \cite{cato04,dono06} \\ 
    Spacings &   [\cite{bap04}, otherwise few so far] \\
    General & Computational tools \cite{koed06}, role in proofs 
\end{tabular}
\end{table}

\subsection{Asymptotic regimes}

\bigskip 
\textbf{ Types of Asymptotics }  \  
The coincidence of ensembles between RMT and statistical theory is
striking, but what can it be \textit{used} for?
The complexity of finite sample size distributions makes the use of
asymptotic approximations appealing, and here an interesting dichotomy
emerges. 
Traditional statistical approximations kept the number of variables
$p$ fixed while letting the sample size $n \to \infty$. This was in
keeping with the needs of the times when 
the number of variables was usually small to moderate.

On the other hand, the nuclear physics models were developed precisely
for settings of high energy levels, and so the number of variables in
the matrix models were large, as seen in the Wigner semi-circle limit.
Interestingly, the many-variables limit of RMT is just what is needed
for modern statistical theories with many variables.

\begin{center}
\begin{tabular}{lcc}
& Stat: $\mathbb{C}$Wishart  &  RMT: Laguerre UE \\  \hline
Density & $\prod_{j=1}^p x_j^{n-p} e^{-x_j} \Delta(x)$ &
$\prod_{j=1}^N x_j^{\alpha} e^{-x_j} \Delta(x)$    \\
\# variables: & $  p$         &     $ N$    \\
Sample size: \qquad  & $n-p$       &     $\alpha$  \\ \hline
\end{tabular}
\end{center}

Comparison of the parameters in the statistics and RMT versions of the
Wishart density in the table above leads to an additional important
remark: in statistics, there is no necessary relationship between
sample size $n$ and number of variables $p$.  We will consider below
limits in which $p/n \to \gamma \in (0, \infty),$ so that $\gamma$
could take any positive value.  In contrast, the most natural
asymptotics in the RMT model would take $N$ large and $\alpha$ fixed.
Thus, from the perspective of orthogonal polynomial theory, the
statistics models lead to somewhat less usual Plancherel-Rotach
asymptotics in which both parameters $N$ \textit{and} $\alpha$ of the
Laguerre polynomials are large.

\bigskip \textbf{Spreading of Sample Eigenvalues  } \ 
To make matters more concrete, we first describe this phenomenon by
example.  Consider $n = 10$ observations on a $p=10$ variable Gaussian
distribution with identity covariance. The sample covariance matrix
follows a Wishart density with $n = p = 10$, and the
\textit{population} eigenvalues $\ell_j(I)$ are all equal to 1.

Nevertheless, there is an extreme spread in the \textit{sample} eigenvalues
$\hat \ell_j = \hat \ell_j(S)$, indeed in a typical sample 
\begin{displaymath}
  (\hat \ell_j) = ( \mathbf{.003}, .036, .095, .16, .30, .51, .78,
  1.12, 1.40, \mathbf{3.07}) 
\end{displaymath}
and the variation is over three orders of magnitude!  Without some
supporting theory, one might be tempted to (erroneously) conclude from
the sample that the population eigenvalues are quite different from
one another.

This spread of sample eigenvalues has long been known, indeed it is an
example of the replusion of eigenvalues induced by the Vandermonde
term in (\ref{eq:jtdens}). It also complicates
the estimation of population covariance matrices -- also a long standing
problem, discussed for example in
\cite{stei77,haff91,yabe94,daka01,ledo04}.

\bigskip \textbf{ The Quarter Circle Law } \ 
\citet{mapa67} gave a systematic
description of the spreading phenomenon: it is the version of the
semi-circle law that applies to sample covariance matrices.
We consider only the special case in which  $A \sim W_p(n, I).$ 
The \textit{empirical distribution function} (or
empirical \textit{spectrum}) counts how many sample eigenvalues fall
below a given value $t$:
\begin{displaymath}
  G_p(t) = p^{-1} \# \{ \hat \ell_j \leq t \}.
\end{displaymath}
The empirical distribution has a
limiting density $g^{MP}$ if sample size $n$ and number of variables $p$ grow
together: $p/n \to \gamma$:
\begin{displaymath}
    g^{MP}(t) = \frac{\sqrt{(b_+-t)(t-b_-)}}{2 \pi \gamma t}, \qquad
 b_{\pm}  = (1 \pm \sqrt{\gamma})^2.
\end{displaymath}

The larger $p$ is relative to $n$, the more spread out is the limiting
density. In particular, with $p = n/4$, one gets the curve supported
in $[\tfrac{1}{4},\tfrac{9}{4}]$.  For $p = n$, the extreme situation
discussed above, the curve covers the full range
from 0 to 4, which corresponds to the huge condition numbers seen in
the sample.

\begin{figure}[h]
  \centering
  \includegraphics[height = 1.5in, width =
  1.5in]{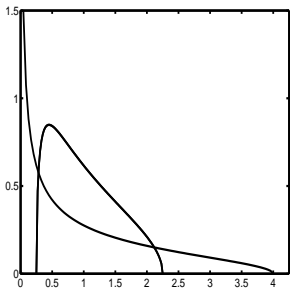} 
  \caption{Mar\v{c}enko-Pastur limit density for $\gamma = \tfrac{1}{4}$
    and $\gamma = 1$.}
  \label{fig:mplaw}
\end{figure}

\section{Largest Eigenvalue Laws}
\label{sec:large-eval}

\textbf{ Hypothesis Test for Largest Eigenvalue } \ 
Suppose that in a sample of $n = 10$ observations from a $p=10$
variate Gaussian distribution $N_{10}(0, \Sigma)$, we see a largest
sample eigenvalue of $4.25$.  Is the observed value consistent with an
identity covariance matrix (with all population eigenvalues = 1), even
though 4.25 lies outside the support interval $[0,4]$ in the
quarter-circle law?

In statistical terms, we are testing a \textit{null hypothesis} of
identity covariance matrix, $H_0: \Sigma = I$ against an
\textit{alternative hypothesis} $H_A: \Sigma \neq I$ that $\Sigma$ has
some more general value.  Normally, of course, one prefers the simpler
model as a description of the data, unless forced by evidence to
conclude otherwise.

One might compare 4.25 to random samples of the largest eigenvalue
from the null hypothesis distribution (three examples yielding
2.91, 3.40 and 3.50); but what is actually needed is an
approximation to the \textup{null hypothesis} distribution of the
largest sample eigenvalue:
\begin{displaymath}
  P\{ \hat \ell_1 > t ~|~ H_0 = W_p(n,I) \}.
\end{displaymath}

\bigskip \textbf{ Tracy-Widom Limits } \ 
Random matrix theory leads to the approximate distribution we need.
In the single Wishart case, assume that $A \sim W_p(n,I),$ either
real or complex, that
$p/n \to \gamma \in (0,\infty)$ and that $\hat \ell_1$ is the largest
eigenvalue in equation (\ref{eq:single-w}).
For the double Wishart case, assume that $A \sim W_p(n_1,I)$ is
independent of $B \sim W_p(n_2,I)$, either real or complex together,
and that $(p/n_1,p/n_2) \to (\gamma_1, \gamma_2) \in (0,1)^2$, and
that $\hat \ell_1$ is the largest generalized eigenvalue in equation
(\ref{eq:double-w}).
With appropriate centering $\mu_{np}$ and scaling $\sigma_{np}$
detailed below, the distribution of the largest eigenvalue
approaches one of the Tracy-Widom $F_\beta$ laws:
\begin{equation}
  \label{eq:Convergence}
     P\{ n \hat \ell_1 \leq \mu_{np} + \sigma_{np} s | H_0 \} \to
   F_\beta(s).
\end{equation}

These laws were first found by Craig Tracy and Harold Widom
\cite{trwi94,trwi96} in the setting of the Gaussian unitary and
orthogonal ensembles, \textit{i.e.} (Hermitian) symmetric Gaussian matrices
with i.i.d.\ entries. There are elegant formulas for the distribution
functions 
$$F_2(s) = \exp \Big( - \int_s^\infty (x-s)^2 q(x) dx \Big), \quad
F_1(s)^2 = F_2(s) \exp \Big( - \int_s^\infty q(x) dx \Big).$$
in terms of the solution $q$ to classical (Painlev\'{e} II) non-linear
second-order differential equation
\begin{displaymath}
  q''= sq + 2 q^3,  \qquad 
  q(s) \sim \mbox{Ai}(s) \ \mbox{as} \ s \rightarrow \infty.
\end{displaymath}
While $q$ and $F_\beta$ are somewhat tricky to compute 
numerically\footnote{At time of writing, for available software in MATLAB 
see \texttt{http://math.arizona.edu/~momar/research.htm} and
\cite{dien05}
in S-PLUS see \texttt{http://www.vitrum.md/andrew/MScWrwck/codes.txt}
and \cite{beja05}. 
Both are based on ideas of \cite{pers02} [see also \cite{edpe05}]
},
from the point of view of applied data
analysis with a software package, it is a special function just like
the normal curve.

As will be seen from the explicit formulas
(\ref{eq:SingComp})-~(\ref{eq:DoubComp1}) below, the scale of
fluctuation $\sigma_{np}/\mu_{np}$ of the largest eigenvalue is
$O(n^{-2/3})$ rather than the $O(n^{-1/2})$ seen in the Gaussian
domain of attraction. This reflects the constraining effect of
eigenvalue repulsion due to the Vandermonde term in (\ref{eq:jtdens}).

The fact that the same limit arises in the single and double Wishart
settings (Laguerre, Jacobi ensembles) is an instance of the
universality discussed in P. Deift's paper \cite{deif06} in this
volume.  In a different direction, one can modify the assumption that
the i.i.d.~entries in the $p \times n$ data matrix $X$ are Gaussian.
Soshnikov \cite{sosh01a} shows that if $n - p = O(p^{1/3})$ and the
matrix entries $X_{ij}$ have sufficiently light (subGaussian) tails,
then the largest eigenvalue continues to have a limiting Tracy-Widom
distribution. The behavior of the largest eigenvalues changes
radically with heavy tailed $X_{ij}$ -- for Cauchy distributed
entries, after scaling by $n^2p^2$,
\cite{sosh06,sofy05} shows a weak form of convergence to a Poisson process.
If the density of the matrix entries behaves like $|x|^{-\mu}$, then
\cite{bbp06} give physical arguments to support a phase transition from
Tracy-Widom to Poisson at $\mu = 4.$

\bigskip \textbf{ Second-order accuracy } \ To demonstrate the
relevance of this limiting result for statistical application, it is
important to investigate its accuracy when the parameters $p$ and $n$
are not so large.  The generic rate of convergence of the left side of
(\ref{eq:Convergence}) to $F_\beta(s)$ is $O(p^{-1/3})$. However,
small modifications in the centering and scaling constants $\mu$ and
$\sigma$, detailed in the four specific cases below, lead to
$O(p^{-2/3})$ errors, which one might call ``second-order accuracy''.
With this improvement, (\ref{eq:Convergence}) takes the form

\begin{equation}
  \label{eq:SecondOrder}
  | P\{ n \hat \ell_1 \leq \mu_{np} + \sigma_{np} s | H_0 \} -
  F_\beta( s) | \leq C  e^{-cs} p^{-2/3}.  
\end{equation}

This higher-order accuracy is reminiscent of that of the central
limit, or normal, approximation to the $t-$test of elementary
statistics for the testing of hypotheses about means, which occurs
when the underlying data has a Gaussian distribution.

\medskip
\textit{Single Wishart, Complex Data.} \   
Convergence in the form (\ref{eq:Convergence}) was first
established by \citet{joha99} as a byproduct of a remarkable analysis
of a random growth model, with
\begin{equation}
  \label{eq:SingComp}
  \mu^o_{np} = ( \sqrt n + \sqrt p )^2, \qquad 
  \sigma^o_{np}  = ( \sqrt n + \sqrt p )
     \biggl( \frac{1}{\sqrt n} + \frac{1}{\sqrt p} \biggr)^{1/3}. 
\end{equation}
The second-order result (\ref{eq:SecondOrder}) is due to
\citet{elka04}, If $\mu_{np}^\prime$ and $\sigma_{np}^\prime$ denote the
quantities in (\ref{eq:SingComp}) with $n$ and $p$ replaced by $n+1/2$
and $p+1/2$, then the centering $\mu_{np}$ is a weighted combination of
$\mu^\prime_{n-1,p}$ and $\mu^\prime_{n,p-1}$ and the scaling $\sigma_{np}$ a similar
combination of $\sigma^\prime_{n-1,p}$ and $\sigma^\prime_{n,p-1}.$

\medskip
\textit{Single Wishart, Real Data.} \ 
Convergence without rates in the form (\ref{eq:Convergence}) to
$F_1(s)$ with centering and scaling as in (\ref{eq:SingComp}) is given
in \cite{john00c}. The assumption that $p/n \to \gamma \in (0,\infty)$
can be weakened to $\min \{n,p\} \to \infty$, as shown by
\citet{elka03} -- this 
extension is of considerable statistical importance since in many
settings $p \gg n$ (see for example \cite{fan06} in these
proceedings).

Analysis along the lines of \cite{john06} suggests that the second
order result (\ref{eq:SecondOrder}) will hold with   
\begin{align}
  \mu_{np} & = \Bigl( \sqrt{\smash[b]{n{  -\hf}}} +
  \sqrt{\smash[b]{p{  -\mathbf{\hf}}}} \Bigr)^2 \label{eq:realmu}\\ 
  \sigma_{np} & = \Bigl( \sqrt{\smash[b]{n{  -\hf}}} +
  \sqrt{\smash[b]{p{  -\hf}}} \Bigr) 
     \biggl( \frac{1}{\sqrt{\smash[b]{n{  -\hf}}}} +
     \frac{1}{\sqrt{\smash[b]{p{  -\hf}}}} \biggr)^{1/3}.  \label{eq:realcov}
\end{align}

\medskip
\textit{Double Wishart, Complex Data.} \ 
Set $\kappa = n_1 + n_2 +1$ and define
\begin{equation}
  \label{eq:DoubComp}
  \sin^2 \Big( \frac{\phi}{2} \Big) = \frac{n_1 + \hf}{\kappa}, \qquad
  \sin^2 \Big( \frac{\gamma}{2} \Big) = \frac{p + \hf}{\kappa}.
\end{equation}
Then
\begin{equation}
  \label{eq:DoubComp1}
  \mu^o_p = \sin^2 \Big( \frac{\phi+\gamma}{2} \Big), \qquad
  (\sigma^o_p)^3 = \frac{ \sin^4(\phi+\gamma)}{4 \kappa^2 \sin \phi
    \sin \gamma}.
\end{equation}
The second-order result (\ref{eq:SecondOrder}) (currently without the
exponential bound, i.e., with $c=0$) is established in \cite{john06}
with $\mu_{np}$ a weighted combination of
$\mu^o_p$ and $\mu^o_{p-1}$ and the scaling $\sigma_{np}$ a similar
combination of $\sigma^o_{p}$ and $\sigma^o_{p-1}.$

\medskip
\textit{Double Wishart, Real Data}. Bound (\ref{eq:SecondOrder}) is shown in
\cite{john06} (again still for $c=0$) with $\mu_{np}$ and
$\sigma_{np}$ given by (\ref{eq:DoubComp1}) with $\kappa = n_1 + n_2 -
1$.

\begin{figure}[h]
  \centering
   \includegraphics[width=2in, height = 2.5in, angle = 0]
 {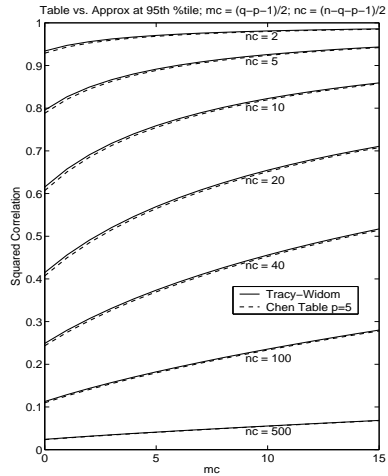}
  \caption{A comparison of the 95th percentile, relevant for
hypothesis tests, from Chen's table (dashed line) and the Tracy-Widom
approximation (solid line). Chen's parameters $m_c, n_c$ are related
to our double Wishart paramaters $n_1, n_2$ by $m_c = (n_1 - p - 1)/2,
n_c = (n_2 - p - 1)/2$.
}
  \label{fig:chenfig}
\end{figure}

\bigskip \textbf{ Approximation vs. Tables for $p = 5$ } \
With second-order correction, Tracy-Widom approximation turns out to
be surprisingly accurate.  William Chen \citep{chen02a, chen03, chen04} has
computed tables of the exact distribution in the double Wishart, real
data, case that cover a wide range of the three parameters $p, n_1$ and
$n_2$, and allow a comparison with the asymptotic approximation.  Even
for $p=5$ variables, the TW approximation is quite good, Figure
\ref{fig:chenfig},  across the
entire range of $n_1$ and $n_2$.

\bigskip \textbf{ A different domain of attraction } \ The Tracy-Widom
laws are quite different from other distributions in the standard
statistical library.  A full probabilistic understanding of their
origin is still awaited (but see \cite{rrv06} for a recent
characterization in terms of the low lying eigenvalues of a random
operator of stochastic diffusion type). Instead, we offer some
incomplete remarks as prelude to the original papers
\cite{trwi94,trwi96}.

Since one is looking at the largest of many eigenvalues, one might be
reminded of extreme value theory, which studies the behavior of the
largest of a collection of variables, which in the simplest case are
independent.  However, extreme value theory exploits the independence to
study the maximum via products: 
$\{ \max_{1 \leq i \leq p} \, l_i \leq t \} = \prod_{i=1}^p I \{
  l_i \leq t \}$
For eigenvalues, however, the Jacobian term, 
or Vandermonde determinant,
\begin{equation} \label{eq:vandermonde}
  V(l) = \prod_{i<j} (l_j - l_i)   = \det [ l_i^{k-1} ]_{1 \leq i,k \leq p},
\end{equation}
changes everything. 
The theory uses the inclusion-exclusion relation:
\begin{displaymath}
  \prod_{i=1}^p  I \{ l_i \leq t \} = \sum_{k=0}^p (-1)^k \binom{p}{k}
\prod_{i=1}^k I\{ l_i > t \}.  
\end{displaymath}
The product structure of the left side, central to extreme value
theory, is discarded in favor of the right side, which
leads to an expression for $P \{ \max_{1 \leq i \leq p} l_i \leq
t \}$ in terms of so-called Fredholm determinants.

For example, it is shown by \citet{trwi98} 
that for complex data 
\begin{displaymath}
  P\{ \max l_i \leq t \}    = \det(I-K_p \chi_{(t,\infty)}),
\end{displaymath}
where $\chi_I$ is the indicator function for interval $I$, and 
$K_p : L_2 \to L_2$ is an operator whose kernel is the two-point
correlation function
\begin{displaymath}
  K_p(x,y) = \sum_{k=1}^p \phi_k(x) \phi_k(y),
\end{displaymath}
written in terms of weighted orthonormal polynomials $\phi_k =
h_k^{-1/2} w^{1/2} p_k$, where the polynomials $p_k$ and weight
functions $w$ are given in Table \ref{tab:polys} for the Gaussian,
Wishart and double Wishart settings respectively.

For real data, \citet{trwi98} show that
\begin{displaymath}
  P\{ \max l_i \leq t \}   = \sqrt{ \det(I-\mathcal{K}_p
    \chi_{(t,\infty)})}, 
\end{displaymath}
where $\mathcal{K}_p$ is now a $2 \times 2$ matrix-valued operator on
$L_2 \otimes L_2$. The corresponding kernel has form
\begin{displaymath}
    \mathcal{K}_p(x,y)  =  
  \begin{pmatrix}
    \tilde K_p   &  -D_2 \tilde K_p  \\
    \epsilon_1 \tilde K_p   &  \tilde K_p^T
  \end{pmatrix},
\end{displaymath}
where $\tilde K_p = K_p + r_1$ and $r_1$ is a rank one kernel
described in the three cases in more detail in 
\cite{afnm00, forr04, john06}. Here $D_2$ and $\epsilon_1$ denote
partial differentiation and integration with respect to second and
first variables respectively.

The expressions are thus somewhat more complicated in the real data case of
primary interest in statistics. However they are amenable to analysis and
approximation using orthogonal polynomial asymptotics near the largest
zero, and to analysis based on the error terms to get the higher order
approximation.

\bigskip \textbf{ Back to the Example } \ 
We asked if
an observed largest eigenvalue of 4.25 was consistent with $H_0:
\Sigma = I$ when $n = p = 10$. The Tracy-Widom
approximation using moments (\ref{eq:realmu})-(\ref{eq:realcov})
yields a 6\% chance of seeing a value more extreme than
4.25 even if ``no structure'' is present.  Against the traditional 5\%
benchmark,  this would not be strong enough evidence to
discount the null hypothesis.

This immediately raises a question about the \textit{power} of the
largest root test, namely evaluation of
\begin{displaymath}
 P \{ \hat \ell_1 > t ~|~ W_p(n,  \Sigma ) \} 
\end{displaymath}
when $\Sigma \neq I$. How different from $1$ does
$\lambda_{\max}(\Sigma)$ need to be before $H_0$ is likely to be rejected?
To this we now turn.

\section{Beyond the Null Hypothesis}
\label{sec:beyond}

From the perspective of multivariate distribution theory, we have, 
in a sense, barely scratched the surface with the classical
RMT ensembles, since they correspond to symmetric situations with no
structure in the population eigenvalues or covariance matrix. 
Basic statistical quantities like power of tests and confidence
intervals, as well as common applications in signal processing,
genetics or finance, call for distributions under 
structured, asymmetric values for the covariance matrix $\Sigma$. 

Statistical theory (pioneered by  Alan James
\cite[e.g.]{jame64}, and summarized in the classic book by Robb
Muirhead \cite{muir82}) gives expressions for the classical
multivariate eigenvalue distributions in more general settings,
typically in terms of hypergeometric functions of matrix
argument.
For example, if $L = \text{diag} (l_i)$ are the eigenvalues of $A \sim
W_p(n,\Sigma)$, then the joint eigenvalue density
\begin{equation*}
  \frac{f_\Sigma(l_1, \ldots, l_p)}{f_I(l_1, \ldots, l_p)} 
     =  |\Sigma|^{-n/2} \exp \{  \tfrac{1}{2} \text{tr} L \} \,
 \text{}_0 \! F_0(- \hf \Sigma^{-1}, L),
\end{equation*}
with
\begin{equation}
  \label{eq:ortho-integral}
  \text{}_0 \! F_0 (S,T)  = \int_{O(p)} \exp \{ \text{tr} ( S U T U^T) \} dU,
\end{equation}
and $dU$ normalized Haar measure, but many other versions occur in the
general theory.  Despite recent major advances in computation by Alan
Edelman and Plamen Koev \cite{koed06,koev06}, and considerable work on the
use of Laplace approximations (see e.g.  \cite{buwo02,buwo05}),
statistical theory would benefit from further serviceable
approximations to these typically rather intractable objects.

\bigskip \textbf{ Persistence of the Tracy-Widom Limit } \
One basic question asks, in the setting of Principal
Components Analysis, for what conditions on the covariance $\Sigma$
does the Tracy-Widom approximation continue to hold, 
\begin{equation}
  \label{eq:twextend}
  P\{\hat \ell_1 \leq \mu_{np}(\Sigma) + \sigma_{np}(\Sigma)s \} \to
  F_\beta(s), 
\end{equation}
perhaps with modified values for centering and scaling to reflect the
value of $\Sigma$?

Fascinating answers are beginning to emerge.
For example, \citet{elka05} establishes that
(\ref{eq:twextend}) holds, along with explicit formulas for
$\mu_{np}(\Sigma)$ and $\sigma_{np}(\Sigma)$,
if enough eigenvalues accumulate near the largest
eigenvalue, or if a small number of eigenvalues are not too isolated,
as we describe below in a specific setting below.

Some of the results are currently restricted to complex data, because
they build in a crucial way on the determinantal representation of the
unitary matrix integral (the complex analog of
(\ref{eq:ortho-integral}))
\begin{equation}
  \label{eq:unitary-integral}
  \int_{U(p)}  \exp \{\text{tr} \Sigma^{-1} U L U^* \} dU
        = c \ \frac{ \det (e^{\pi_j l_k}) }{V(\pi) V(l)} 
\end{equation}
known as the Harish-Chandra-Itzykson-Zuber formula \cite{hc57,itzu80},
see also \cite{grri89}.  Here the eigenvalues of $\Sigma^{-1}$ are
given by $\text{diag} (\pi_j)$ and $V(l)$
is the Vandermonde determinant (\ref{eq:vandermonde}).  While it is
thought unlikely that there are direct analogs of
(\ref{eq:unitary-integral}), we very much need extensions of the
distributional results to real data: there are some results in the
physics literature \cite{brhi03}, but any statistical consequences are
still unclear.

\bigskip \textbf{ Finite rank perturbations.} \ We focus on a simple
concrete model, and describe a phase transition phenomenon.  Assume
that 
\begin{equation}
  \label{eq:spike}
\Sigma = \text{diag}(\ell_1, \ldots, \ell_M, \sigma^2_e, \ldots,
\sigma^2_e),  
\end{equation}
so that a fixed number $M$ of population eigenvalues are greater than
the base level $\sigma^2_e$, while both dimensions $p$ and $n$
increase in constant ratio $ p/n \to \gamma \in (0,\infty).$

First some heuristics: if all population eigenvalues are equal, then the
largest sample eigenvalue $\hat \ell_1$ has $n^{-2/3}$ fluctuations
around the upper limit of the support of the Mar\v{c}enko-Pastur quarter
circle law, the fluctuations being described by the Tracy-Widom law. 
For simplicity, consider $M=1$ and $\sigma^2_e = 1$. 
If $\ell_1$ is large and so very clearly separated from the bulk
distribution, then one expects Gaussian fluctuations of order
$n^{-1/2}$, and this is confirmed by standard perturbation analysis.
 
\citet{bbap05} describe, for \textit{complex} data, a `phase
transition' that occurs between these two extremes.  If $\ell_1 \leq 1
+ \sqrt \gamma$, then
\begin{displaymath}
  n^{2/3}(\hat \ell_1 - \mu)/\sigma  \Rightarrow
  \begin{cases}
    F_2 \qquad & \ell_1 < 1 + \sqrt \gamma \\
    \tilde F_2 & \ell_1 = 1 + \sqrt \gamma 
  \end{cases}
\end{displaymath}
where, from (\ref{eq:SingComp}), we may set
\begin{displaymath}
  \mu = (1 + \sqrt \gamma)^2, \qquad
  \sigma = (1 + \sqrt \gamma) (1 + \sqrt{ \gamma^{-1}} )^{1/3},
\end{displaymath}
and $\tilde F_2$ is related to $F_2$ as described in
\citet{bbap05}. 
On the other hand, if $\ell_1 > 1 + \sqrt \gamma,$
\begin{gather}
  n^{1/2} ( \hat \ell_1 - \mu(\ell_1))/\sigma(\ell_1)
  \Rightarrow N(0,1), \notag \\
\intertext{with}
  \mu(\ell_1) = \ell_1 \Bigl(1 + \frac{\gamma}{\ell_1 - 1} \Bigr),
\qquad
\sigma^2(\ell_1) =  \ell_1^2 \Bigl(1 - \frac{\gamma}{(\ell_1 -
  1)^2} \Bigr).  \label{eq:moments}
\end{gather}

Thus, below the phase transition the distribution of
$\hat \ell_1$ is unchanged, Tracy-Widom, regardless of the value of
$\ell_1 $. 
As $\ell_1$ increases through $1 + \sqrt \gamma$, the law of $\hat \ell_1$
jumps to Gaussian and the mean increases with $\ell_1$, but 
is biased low, $\mu(\ell_1) < \ell_1,$ while the variance 
$\sigma^2(\ell_1)$ is lower than 
its value, $\ell_1^2$, in the limit with $p$ fixed.

A key feature is that the phase transition point $1 + \sqrt \gamma$,
located at the zero of $\sigma(\ell_1$), is buried deep inside the
bulk, whose upper limit is $(1 + \sqrt \gamma)^2$. A good heuristic
explanation for this location is still lacking, though see
\citet{elka05}.

Further results on almost sure and Gaussian limits for both real and
complex data, and under weaker distributional assumptions have been
obtained by \citet{paul04} and \citet{basi06}.

\bigskip \textbf{ A recent example. }  \ \citet{hard06} illustrates
simply this phase transition phenomenon in a setting from economics
and finance.  In a way this is a negative example for PCA; but
statistical theory is as concerned with describing the limits of
techniques as their successes.

Factor analysis models, of recently renewed
interest in economics, attempt to ``explain'' the prices or returns of
a portfolio of securities in terms of a small number of common
``factors'' combined with security-specific noise terms.
It has been further postulated that one could estimate the number and
sizes of these factors using PCA.  In a 1989 paper that is widely
cited and taught in economics and finance, Brown
\cite{brow89} gave a realistic simulation example that challenged this
view, in a way that remained incompletely understood until recently.

\begin{figure}[h]
\parbox{.49\textwidth}{
  \centering
  \includegraphics[width =
  .48\textwidth]{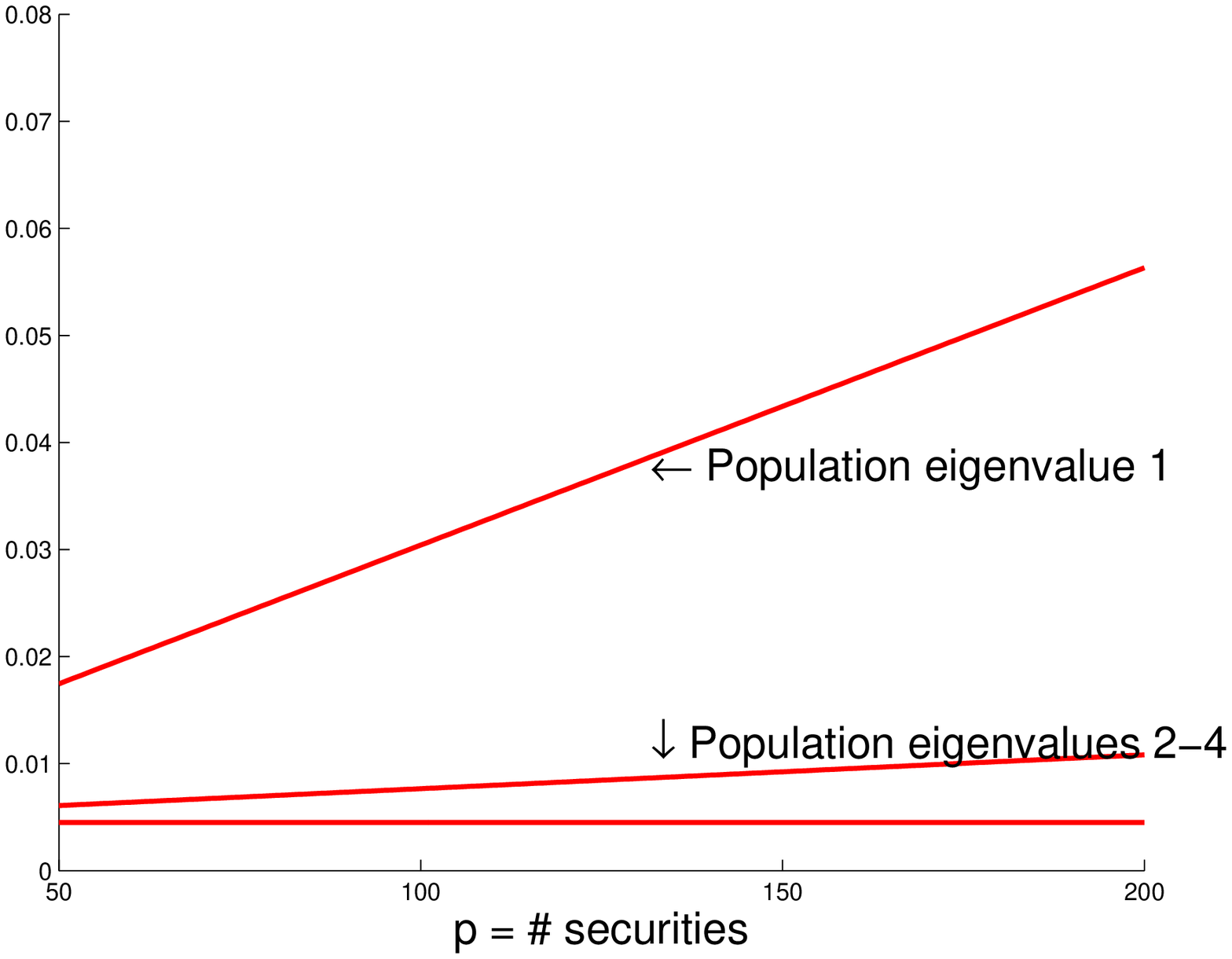} }
\parbox{.49\textwidth}{
  \centering
  \includegraphics[width =
  .48\textwidth]{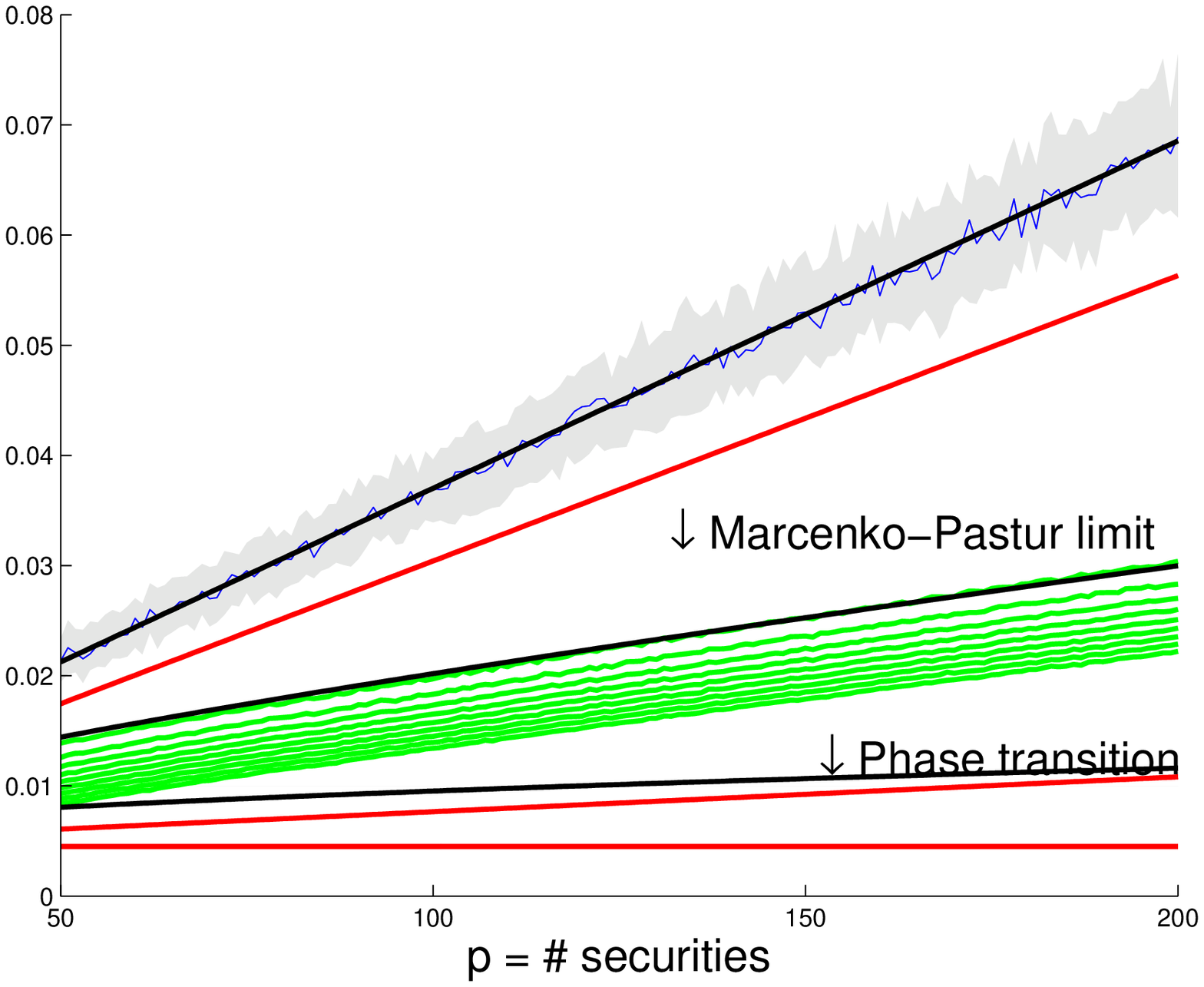} }
\caption{Population and sample eigenvalues for a four factor model
  (\ref{eq:brownmodel}) with $\beta = 0.6, \sigma_b = .4, \sigma_f =
  .01257, \sigma_e = .0671.$ [Brown \& Harding use $\beta = 1,
  \sigma_b = .1$; the values are modified here for legibility of the
  plot.]  (a) Left panel: Population eigenvalues according to
  (\ref{eq:brown-evals}) (b) Right panel: The top sample eigenvalue in
  replications spreads about a sample average line which tracks the
  solid line given by (\ref{eq:moments}), in particular overestimating
  the population value $\ell_1$. The next nine sample eigenvalues fall
  at or below the Mar\v{c}enko-Pastur upper limit, swamping the next
  three population eigenvalues.  }
  \label{fig:harding}
\end{figure}

Brown's example postulated four independent factors, with the parameters
of the model calibrated to historical data from the New York Stock Exchange.
The return in period $t$ of security $k$ is assumed to be given by
\begin{equation} \label{eq:brownmodel}
  R_{kt} =  \Sigma_{\nu =1}^4 b_{k\nu } f_{\nu t} + e_{kt}; 
\qquad k  = 1, \ldots, p; \ \ t = 1, \ldots, T,
\end{equation}
where it is assumed that $b_{k\nu } \sim N(\beta, \sigma_b^2),
f_{\nu t} \sim N(0, \sigma_f^2)$ and 
$  e_{\nu t} \sim N(0, \sigma_e^2)$, all independently of one another.
The population covariance matrix has the form 
(\ref{eq:spike}) with
$M=4$ and
\begin{equation}
  \label{eq:brown-evals}
  \ell_j = p \sigma^2_f (\sigma^2_b + 4 \beta \delta_{j1}) +
\sigma^2_e,  \qquad \qquad j = 1, \ldots, 4.
\end{equation} 
Here $\delta_{j1}$ is the Kronecker delta, equal to $1$ for $j=1$ and
$0$ otherwise.
Figure \ref{fig:harding}(a) plots the population eigenvalues $\ell_1$
(the dominant `market' factor), the common value $\ell_2 = \ell_3 =
\ell_4$ and the base value $\ell_5 = \sigma^2_e$ against $p$, the
number of securities in the portfolio.  One might expect to be able to
recover an estimate of $\ell_2$ from empirical data, but this turns
out to be impossible for $p \in [50, 200]$ when $T = 80$ as shown in
Figure \ref{fig:harding}(b).  First, the range of observed values of
the top or market eigenvalue is biased upward from the true top
eigenvalue.  In addition, there are many sample eigenvalues above the
anticipated value for $\ell_2$.  

Harding shows that one can directly
apply the (real version) of the phase transition results previously
discussed to fully explain Brown's results.
Indeed, the inability to identify factors is because they fall on the
wrong side of the phase transition $\sigma^2_e(1 + \sqrt{p/T})$, and
so we can not expect the observed eigenvalue estimates to exceed the
Mar\v{c}enko-Pastur upper bound $\sigma^2_e(1 + \sqrt{p/T})^2$. Finally,
the bias between the observed and true values of the top eigenvalue is
also accurately predicted by the random matrix formulas
(\ref{eq:moments}).

\section{Estimating Eigenvectors}
\label{sec:eigenvector-est}

Most of the literature at the intersection of random matrix theory and
statistics is focused on eigenvalues. We close with a few remarks on
the estimation of eigen\textit{vectors}. Of course, the question is
only meaningful in non-symmetric settings when the covariance matrix
$\Sigma$ is not proportional to $I$. We again assume that $S \sim
W_p(n,\Sigma)$ and now focus attention on covariance
models which are a finite-rank perturbation of the identity\footnote{
The situation is different in \textit{functional} Principal Components
Analysis, where smoothness of the observed data (functions) leads to
covariance matrices with smoothly decaying eigenvalues. For entries
into this literature, see for example \cite{dpr82, bosq00, hmw06}}:
\begin{equation}
  \label{eq:spikemx}
  \Sigma = \sigma^2 I + \sum_{\nu=1}^M \lambda_\nu \theta_\nu \theta_\nu^T,
\end{equation}
with $\lambda_1 \geq \ldots \geq \lambda_M > 0$ and $\{ \theta_\nu \}$
orthonormal. We ask how well can the population eigenvectors
$\theta_\nu$ be estimated when both $p$ and $n$ are large.

First some remarks on how model (\ref{eq:spikemx}) can arise from an
\textit{orthogonal factor} or \textit{variance components} model for
the data.  Assume that the $p-$dimensional observations $X_i, i= 1,
\ldots ,n$ have the form
\begin{equation}
  \label{eq:var-comp}
  X_i = \mu + 
\sum_{\nu=1}^M \sqrt{\lambda_\nu} v_{\nu i} \theta_\nu + \sigma Z_i,
\notag
\end{equation}
where $\{v_{\nu i}: 1\leq \nu \leq M\}$ are i.i.d. $N(0,1),$
independently of $Z_i \sim N_p(0,I_p)$, for all $i$.
If we further assume, for convenience, that $\mu = 0$, then with the
sample covariance $S$ defined as in Section \ref{sec:Background}, then
$S \sim W_p(n,\Sigma)$.
If we express $X_i, \theta_\nu$ and $Z_i$ in
(\ref{eq:var-comp}) in terms of co-ordinates in a
suitable basis $\{ e_k, k = 1, \ldots, p \}$ and write $f_{\nu i} =
\sqrt{ \lambda_\nu} v_{\nu i}$  we obtain
\begin{displaymath}
  X_{ki} = \sum_{\nu =1}^M \theta_{k \nu} f_{\nu i} +\sigma Z_{ki},
\end{displaymath}
in which $\theta_{k \nu}$ is viewed as the factor loading of the $k$th
variable on the $\nu$th factor, and $f_{\nu i}$ is the factor score of
the $\nu$th factor for the $i$th individual.  As we have seen in
(\ref{eq:brownmodel}) in the previous section, in economics $X_{ki}$
may represent the return on the $k$th security in time period $i$.

Assume that $\lambda_1 > \ldots > \lambda_M > 0$. Let $\hat
\theta_\nu$ denote the normalized sample eigenvectors of $S$ (denoted
$\widehat{\mathbf{v}}_\nu$ in Section \ref{sec:PCA}) associated with
the $M$ largest sample eigenvalues. 
In classical asymptotics, with $n$ large and $p$ fixed, there is a
well understood Gaussian limit theory:
\begin{equation}
  \label{eq:classical-limit}
  \sqrt{n} (\hat \theta_\nu - \theta_\nu) \to N_p(0, \Gamma_\nu)
\end{equation}
where $\Gamma_\nu$ is given, for example, in \cite{ande63,ande84}.

The situation is radically different when $p/n \to \gamma > 0$ 
-- indeed, ordinary PCA is necessarily inconsistent:
\begin{equation*}
 \langle \hat \theta_\nu, \theta_\nu \rangle \rightarrow 
  \begin{cases}
    0 & \lambda_\nu \in [0,\sqrt{\gamma}] \\
    \frac{1 - \gamma/\lambda_\nu^2}{1+ \gamma/\lambda_\nu} &
    \lambda_\nu > \sqrt \gamma 
  \end{cases},
\end{equation*}
For signal strengths $\lambda$ below the phase transition just
discussed, nothing can be estimated -- the estimate is asymptotically
orthogonal to the truth. The angle decreases as $\lambda_\nu$ grows, but
is never exactly consistent.

This result has emerged in several literatures, starting in the
learning theory/statistical physics community, with non-rigorous
arguments based on the replica method \cite{rebb96, hora04}, where 
this phenomenon has been termed ``retarded learning'' 
\cite{bimi94,wana94}. 
More recently, rigorous results
have been obtained \cite{jolu04, paul04, onat06b}.

To obtain consistent estimates, further assumptions are needed. One
plausible situation is that in which there exists a basis $\{ e_k
\}_{k=1:p}$ in which it is believed that the vectors $\theta_\nu$ have a
sparse representation. In microarray genetics, for example $X_{ki}$
might be the expression of gene $k$ in the $i$th patient, and it may
be believed that (in the standard basis) each factor $\nu$ is related
to only a small number of genes \cite{lcwm06}. In EEG studies of the
heart, the beat-to-beat cycle might be expressed in a wavelet basis,
in which the components of variation $\theta_\nu$ may well be
sparsely represented \cite{jolu04}.

We briefly describe results in the sparse setting of work in progress
by D. Paul, and by Paul and the author. For simplicity only, we
specialize to $M=1$.  The error of estimation, or loss, of $\hat
\theta$ is measured on unit vectors by
\begin{displaymath}
  L(\hat \theta,\theta) = \| \hat \theta - \text{sign}(\langle \hat
  \theta, \theta \rangle) \theta \|^2 = 4 \sin^2 \tfrac{1}{2} \angle (\hat
  \theta, \theta).
\end{displaymath}
If $\hat \theta$ is now the ordinary PCA estimate of $\theta$, and if
$p/n \to \gamma > 0$, then to first order,
\begin{displaymath}
  E L(\hat \theta, \theta) = \frac{p}{n h(\lambda)} (1 + o(1)), \qquad
  \qquad  
  h(\lambda) = \frac{\lambda^2}{1+\lambda},
\end{displaymath}
from which it is natural to define the ``per-variable'' noise level
$\tau_n = 1/\sqrt{n h(\lambda)}$.

As is common in non-parametric estimation theory, we use $\ell_q$
norm, $q < 2$, as a measure of sparsity: with $\| \theta \|_q^q =
\sum_k |\theta_k|^2$, define $\Theta_q(C) = \{ \theta \in S^{p-1} ~:~ \|
\theta \|_q \leq C \}$. 
Paul proposes a two-step procedure for selecting a reduced subset of
variables on which to perform PCA, resulting in an estimator
$\hat{\theta}^P$ for which
\begin{equation} \label{eq:sparse-est}
  \sup_{ \theta \in \Theta_q(C)} E L(\hat \theta^P, \theta) 
      \leq K(C) \log p \cdot m_n \tau_n^2.
\end{equation}
Here $m_n$ is an effective dimension parameter, equal to
$(C^2/(\tau^2 \log p))^{q/2}$ in the ``sparse'' case when this is
smaller than $c_1 p$, and equal to $p$ in the contrary ``dense''
case. 
Lower bounds are obtained that show that this estimation error is
optimal, in a minimax sense, up to factors that are at most
logarithmic in $p$.

Bounds such as (\ref{eq:sparse-est}) are reminiscent of those for
estimation of sparse \textit{mean} sequences in white Gaussian noise
\cite{dojo94a,bima99,josi02a}. 
An observation due to Paul provides a link between
eigenvector estimation and the estimation of means. Again with $M=1$
for simplicity, let $\hat \theta$ be the ordinary PCA estimate of
$\theta$. Write $\hat C = \langle \hat \theta, \theta \rangle$ and
$\hat \theta^\perp = \hat \theta -  C \theta$. 
Then, with $\hat S^2 = 1 - \hat C^2$, in the decomposition
\begin{displaymath}
  \hat \theta = \hat C \theta + \hat S U, \qquad \qquad 
  U = \hat \theta^\perp / \| \hat \theta^\perp \|
\end{displaymath}
it happens that $U$ is uniformly distributed on a copy of $S^{p-2}$,
independently of $\hat S$. 

It is a classical remark that a
high-dimensional isotropic Gaussian vector is essentially concentrated
uniformly on a sphere. We may reverse this remark by starting with a
uniform distribution on a sphere, and introducing an ultimately
inconsequential randomization with $R^2 \sim \chi_{p-1}^2/p$ and $z_1
\sim N(0,1/p)$ with the result that $z = R U + z_1 \theta$ has an
$N_p(0,I)$ distribution. 
This leads to a signal-in-Gaussian-noise representation
\begin{displaymath}
  Y = \hat C \theta + \tau^2 z,  \qquad \qquad 
   \tau^2 = 1/(2 n h(\hat \lambda)),
\end{displaymath}
Work is in progress to use this connection
to improve the extant estimation results for eigenvectors.

\section{Coda}
\label{sec:coda}

One may expect a continuing fruitful influence of developments in
random matrix theory on high dimensional statistical theory, and
perhaps even some flow of ideas in the opposite direction. A snapshot
of current trends may be obtained from  
\url{http://www.samsi.info/workshops/2006ranmat-opening200609.shtml}, being
the presentations from the Opening Workshop of a semester devoted to
High Dimensional Inference and Random Matrices at the NSF Statistics
and Applied Mathematics Institute in Fall 2006.


\begin{thebibliography}{98}
\providecommand{\natexlab}[1]{#1}
\providecommand{\url}[1]{\texttt{#1}}
\expandafter\ifx\csname urlstyle\endcsname\relax
  \providecommand{\doi}[1]{doi: #1}\else
  \providecommand{\doi}{doi: \begingroup \urlstyle{rm}\Url}\fi

\bibitem[Adler et~al.(2000)Adler, Forrester, Nagao, and van Moerbeke]{afnm00}
M.~Adler, P.~J. Forrester, T.~Nagao, and P.~van Moerbeke.
\newblock Classical skew orthogonal polynomials and random matrices.
\newblock \emph{Journal of Statistical Physics}, 99\penalty0 (1/2):\penalty0
  141--170, 2000.

\bibitem[Anderson(1963)]{ande63}
T.~W. Anderson.
\newblock Asymptotic theory for principal component analysis.
\newblock \emph{Annals of Mathematical Statistics}, 34:\penalty0 122--148,
  1963.

\bibitem[Anderson(1984)]{ande84}
T.~W. Anderson.
\newblock \emph{An Introduction to Multivariate Statistical Analysis, 2nd ed.}
\newblock Wiley, 1984.

\bibitem[Andersson et~al.(1983)Andersson, Br{\o}ns, and Jensen]{abj83}
Steen~A. Andersson, Hans~K. Br{\o}ns, and S{\o}ren~Tolver Jensen.
\newblock Distribution of eigenvalues in multivariate statistical analysis.
\newblock \emph{Ann. Statist.}, 11\penalty0 (2):\penalty0 392--415, 1983.

\bibitem[Bai(1999)]{bai99}
Z.~D. Bai.
\newblock Methodologies in spectral analysis of large dimensional random
  matrices, a review.
\newblock \emph{Statistica Sinica}, 9:\penalty0 611--677, 1999.

\bibitem[Baik and Silverstein(2006)]{basi06}
Jinho Baik and Jack~W. Silverstein.
\newblock Eigenvalues of large sample covariance matrices of spiked population
  models.
\newblock \emph{Journal of Multivariate Analysis}, 97:\penalty0 1382--1408,
  2006.

\bibitem[Baik et~al.(2005)Baik, Ben~Arous, and P{\'e}ch{\'e}]{bbap05}
Jinho Baik, G{\'e}rard Ben~Arous, and Sandrine P{\'e}ch{\'e}.
\newblock Phase transition of the largest eigenvalue for nonnull complex sample
  covariance matrices.
\newblock \emph{Ann. Probab.}, 33\penalty0 (5):\penalty0 1643--1697, 2005.

\bibitem[Barnett and Preisendorfer(1987)]{bapr87}
T.~P. Barnett and R.~Preisendorfer.
\newblock Origins and levels of monthly and seasonal forecast skill for
  {U}nited {S}tates surface air temperatures determined by canonical
  correlation analysis.
\newblock \emph{Monthly Weather Review}, 115:\penalty0 1825--1850, 1987.

\bibitem[Bejan(2005)]{beja05}
A.~Bejan.
\newblock Largest eigenvalues and sample covariance matrices. {T}racy-{W}idom
  and {P}ainlev¶e {II}: computational aspects and realization in {S}-{P}lus
  with applications.
\newblock \url{http://www.vitrum.md/andrew/TWinSplus.pdf}, 2005.

\bibitem[Ben~Arous and P{\'e}ch{\'e}(2005)]{bap04}
G.~Ben~Arous and S.~P{\'e}ch{\'e}.
\newblock Universality of local eigenvalue statistics for some sample
  covariance matrices.
\newblock \emph{Comm. Pure Appl. Math.}, 58\penalty0 (10):\penalty0 1316--1357,
  2005.

\bibitem[Biehl and Mietzner(1994)]{bimi94}
M.~Biehl and A.~Mietzner.
\newblock Statistical mechanics of unsupervised structure recognition.
\newblock \emph{Journal of Physics A: Mathematical and General}, 27\penalty0
  (6):\penalty0 1885--1897, 1994.

\bibitem[Birg{\'e} and Massart(2001)]{bima99}
Lucien Birg{\'e} and Pascal Massart.
\newblock Gaussian model selection.
\newblock \emph{Journal of European Mathematical Society}, 3:\penalty0
  203--268, 2001.

\bibitem[Biroli et~al.(2006)Biroli, Bouchaud, and Potters]{bbp06}
Giulio Biroli, Jean-Philippe Bouchaud, and Marc Potters.
\newblock On the top eigenvalue of heavy-tailed random matrices, 2006.
\newblock arXiv:cond-mat/0609070.

\bibitem[Bosq(2000)]{bosq00}
D.~Bosq.
\newblock \emph{Linear processes in function spaces}, volume 149 of
  \emph{Lecture Notes in Statistics}.
\newblock Springer-Verlag, New York, 2000.

\bibitem[Bouchaud and Potters(2003)]{bopo03}
Jean-Philippe Bouchaud and Marc Potters.
\newblock \emph{Theory of Financial Risk and Derivative Pricing: From
  Statistical Physics to Risk Management}.
\newblock Cambridge University Press, 2003.

\bibitem[Box(1979)]{box79}
G.~E.~P. Box.
\newblock Robustness in the strategy of scientific model building.
\newblock In R.~L. Launer and G.~N. Wilkinson, editors, \emph{Robustness in
  Statistics}. Academic Press: New York, 1979.

\bibitem[Br\'{e}zin and Hikami(2003)]{brhi03}
E.~Br\'{e}zin and S.~Hikami.
\newblock New correlation functions for random matrices and integrals over
  supergroups.
\newblock \emph{Journal of Physics A: Mathematical and General}, 36\penalty0
  (3):\penalty0 711--751, 2003.

\bibitem[Brown(1989)]{brow89}
Stephen~J. Brown.
\newblock The number of factors in security returns.
\newblock \emph{The Journal of Finance}, XLIV\penalty0 (5):\penalty0
  1247--1261, 1989.

\bibitem[Butler and Wood(2002)]{buwo02}
Ronald~W. Butler and Andrew T.~A. Wood.
\newblock Laplace approximations for hypergeometric functions with matrix
  argument.
\newblock \emph{Ann. Statist.}, 30\penalty0 (4):\penalty0 1155--1177, 2002.

\bibitem[Butler and Wood(2005)]{buwo05}
Ronald~W. Butler and Andrew T.~A. Wood.
\newblock Laplace approximations to hypergeometric functions of two matrix
  arguments.
\newblock \emph{J. Multivariate Anal.}, 94\penalty0 (1):\penalty0 1--18, 2005.

\bibitem[Cand\`{e}s and Tao(2004)]{cato04}
Emmanuel Cand\`{e}s and Terence Tao.
\newblock {Near Optimal Signal Recovery From Random Projections: Universal
  Encoding Strategies?}, 2004.
\newblock arXiv:math.CA/0410542.

\bibitem[Cavalli-Sforza(2000)]{cs00}
L.~L. Cavalli-Sforza.
\newblock \emph{Genes, peoples, and languages}.
\newblock North Point Press, 2000.

\bibitem[Cavalli-Sforza et~al.(1994)Cavalli-Sforza, Menozzi, and
  Piazza]{csmp94}
L.~Luca Cavalli-Sforza, Paolo Menozzi, and Alberto Piazza.
\newblock \emph{The history and geography of human genes}.
\newblock Princeton University Press, 1994.

\bibitem[Chen(2002)]{chen02a}
William~R. Chen.
\newblock Some new tables of the largest root of a matrix in multivariate
  analysis: A computer approach from 2 to 6, 2002.
\newblock Presented at the 2002 American Statistical Association.

\bibitem[Chen(2003)]{chen03}
William~W. Chen.
\newblock Table for upper percentage points of the largest root of a
  determinantal equation with five roots.
\newblock \emph{InterStat}, \penalty0 (5), February 2003.
\newblock URL \url{interstat.statjournals.net}.

\bibitem[Chen(2004)]{chen04}
William~W. Chen.
\newblock The new table for upper percentage points of the largest root of a
  determinantal equation with seven roots.
\newblock \emph{InterStat}, \penalty0 (1), September 2004.
\newblock URL \url{interstat.statjournals.net}.

\bibitem[Daniels and Kass(2001)]{daka01}
Michael~J. Daniels and Robert~E. Kass.
\newblock Shrinkage estimators for covariance matrices.
\newblock \emph{Biometrics}, 57\penalty0 (4):\penalty0 1173--1184, 2001.

\bibitem[Dauxois et~al.(1982)Dauxois, Pousse, and Romain]{dpr82}
J.~Dauxois, A.~Pousse, and Y.~Romain.
\newblock Asymptotic theory for the principal component analysis of a vector
  random function: some applications to statistical inference.
\newblock \emph{J. Multivariate Anal.}, 12\penalty0 (1):\penalty0 136--154,
  1982.

\bibitem[Deift(2007)]{deif06}
P.~Deift.
\newblock Universality for mathematical and physical systems.
\newblock In \emph{Proceedings of the International Congress of
  Mathematicians}, volume~I, pages XXX--XXX. EMS Publishing House, Z\"{u}rich,
  2007.

\bibitem[Diaconis(2003)]{diac03}
Persi Diaconis.
\newblock Patterns in eigenvalues: the 70th {J}osiah {W}illard {G}ibbs lecture.
\newblock \emph{Bull. Amer. Math. Soc. (N.S.)}, 40\penalty0 (2):\penalty0
  155--178 (electronic), 2003.

\bibitem[Dieng()]{dien05}
Momar Dieng.
\newblock {Distribution Functions for Edge Eigenvalues in Orthogonal and
  Symplectic Ensembles: Painlev\'e Representations II}.
\newblock arXiv:math.PR/0506586.

\bibitem[Donoho and Johnstone(1994)]{dojo94a}
D.~L. Donoho and I.~M. Johnstone.
\newblock Ideal spatial adaptation via wavelet shrinkage.
\newblock \emph{Biometrika}, 81:\penalty0 425--455, 1994.

\bibitem[Donoho(2006)]{dono06}
David~L. Donoho.
\newblock For most large underdetermined systems of linear equations the
  minimal {$l\sb 1$}-norm solution is also the sparsest solution.
\newblock \emph{Comm. Pure Appl. Math.}, 59\penalty0 (6):\penalty0 797--829,
  2006.

\bibitem[Dyson(1962)]{dyso62}
Freeman~J. Dyson.
\newblock The threefold way. {A}lgebraic structure of symmetry groups and
  ensembles in quantum mechanics.
\newblock \emph{Journal of Mathematical Physics}, 3\penalty0 (6):\penalty0
  1199--1215, 1962.

\bibitem[Edelman and Persson()]{edpe05}
Alan Edelman and Per-Olof Persson.
\newblock {Numerical Methods for Eigenvalue Distributions of Random Matrices}.
\newblock arXiv:math-ph/0501068.

\bibitem[Edelman and Rao(2005)]{edra05}
Alan Edelman and N.~Raj Rao.
\newblock Random matrix theory.
\newblock \emph{Acta Numer.}, 14:\penalty0 233--297, 2005.

\bibitem[El~Karoui(2003)]{elka03}
Noureddine El~Karoui.
\newblock {On the largest eigenvalue of {W}ishart matrices with identity
  covariance when $n, p$ and $p/n$ tend to infinity}, 2003.
\newblock arXiv:math.ST/0309355.

\bibitem[El~Karoui(2004)]{elka04}
Noureddine El~Karoui.
\newblock {An asymptotic Berry-Esseen result for the largest eigenvalue of
  complex white {W}ishart matrices}, 2004.
\newblock arXiv:math.PR/0409610.

\bibitem[El~Karoui(2005)]{elka05}
Noureddine El~Karoui.
\newblock Tracy-{W}idom limit for the largest eigenvalue of a large class of
  complex {W}ishart matrices, 2005.
\newblock arXiv:math.PR/0503109, To appear in {\it Annals of Probability}.

\bibitem[Fan and Li(2006)]{fan06}
J.~Fan and R.~Li.
\newblock {Statistical Challenges with High Dimensionality: Feature Selection
  in Knowledge Discovery}.
\newblock In \emph{Proceedings of the International Congress of
  Mathematicians}, volume III, pages 595--622. EMS Publishing House,
  Z\"{u}rich, 2006.

\bibitem[Fisher(1939)]{fish39}
R.~A. Fisher.
\newblock The sampling distribution of some statistics obtained from non-linear
  equations.
\newblock \emph{Annals of Eugenics}, 9:\penalty0 238--249, 1939.

\bibitem[Forrester(2004)]{forr04}
P.~J. Forrester.
\newblock Log-gases and {R}andom matrices.
\newblock URL \url{http://www.ms.unimelb.edu.au/~matpjf/matpjf.html}.
\newblock Book manuscript, 2004.

\bibitem[Fox and Kahn(1964)]{foka64}
David Fox and Peter~B. Kahn.
\newblock Higher order spacing distributions for a class of unitary ensembles.
\newblock \emph{Phys. Rev.}, 134\penalty0 (5B):\penalty0 B1151--B1155, Jun
  1964.

\bibitem[Friman et~al.(2001)Friman, Cedefamn, Lundberg, Borga, and
  Knutsson]{frim01}
Ola Friman, J.~Cedefamn, P.~Lundberg, M.~Borga, and H.~Knutsson.
\newblock Detection of neural activity in functional {MRI} using canonical
  correlation analysis.
\newblock \emph{Magnetic Resonance in Medicine}, 45:\penalty0 323--330, 2001.

\bibitem[Girshick(1939)]{girs39}
M.~A. Girshick.
\newblock On the sampling theory of roots of determinantal equations.
\newblock \emph{Annals of Mathematical Statistics}, 10:\penalty0 203--224,
  1939.

\bibitem[Gross and Richards(1989)]{grri89}
Kenneth~I. Gross and Donald St.~P. Richards.
\newblock Total positivity, spherical series, and hypergeometric functions of
  matrix argument.
\newblock \emph{J. Approx. Theory}, 59\penalty0 (2):\penalty0 224--246, 1989.

\bibitem[Haff(1991)]{haff91}
L.~R. Haff.
\newblock The variational form of certain {B}ayes estimators.
\newblock \emph{Annals of Statistics}, 19:\penalty0 1163--1190, 1991.

\bibitem[Hall et~al.(2006)Hall, M\"uller, and Wang]{hmw06}
Peter Hall, Hans-Georg M\"uller, and Jane-Ling Wang.
\newblock Properties of principal component methods for functional and
  longitudinal data analysis.
\newblock \emph{Annals of Statistics}, 34\penalty0 (3):\penalty0 1493--1517,
  2006.

\bibitem[Harding(2006)]{hard06}
Matthew~C. Harding.
\newblock Explaining the single factor bias of arbitrage pricing models in
  finite samples.
\newblock URL \url{http://www.mit.edu/~mharding/}.
\newblock Dept. of Economics, MIT, 2006.

\bibitem[Harish-Chandra(1957)]{hc57}
Harish-Chandra.
\newblock Differential operators on a semisimple {L}ie algebra.
\newblock \emph{American Journal of Mathematics}, 79\penalty0 (1):\penalty0
  87--120, 1957.

\bibitem[Hotelling(1933)]{hote33}
H.~Hotelling.
\newblock Analysis of a complex of statistical variables into principal
  components.
\newblock \emph{J. Educational Psychology}, 24:\penalty0 417--441,498--520,
  1933.

\bibitem[Hotelling(1936)]{hote36}
H.~Hotelling.
\newblock Relations between two sets of variates.
\newblock \emph{Biometrika}, 28:\penalty0 321--377, 1936.

\bibitem[Hoyle and Rattray(2004)]{hora04}
D.~C. Hoyle and M.~Rattray.
\newblock Principal-component-analysis eigenvalue spectra from data with
  symmetry breaking structure.
\newblock \emph{Physical Review E}, 69\penalty0 (026124), 2004.

\bibitem[Hsu(1939)]{hsu39}
P.~L. Hsu.
\newblock On the distribution of roots of certain determinantal equations.
\newblock \emph{Annals of Eugenics}, 9:\penalty0 250--258, 1939.

\bibitem[Itzykson and Zuber(1980)]{itzu80}
C.~Itzykson and J.-B. Zuber.
\newblock The planar approximation. {II}.
\newblock \emph{J. Mathematical Physics}, 21\penalty0 (3):\penalty0 411--421,
  1980.

\bibitem[James(1964)]{jame64}
A.~T. James.
\newblock Distributions of matrix variates and latent roots derived from normal
  samples.
\newblock \emph{Annals of Mathematical Statistics}, 35:\penalty0 475--501,
  1964.

\bibitem[Johansson(2000)]{joha99}
Kurt Johansson.
\newblock Shape fluctuations and random matrices.
\newblock \emph{Communications in Mathematical Physics}, 209:\penalty0
  437--476, 2000.

\bibitem[Johnson and Wichern(2002)]{jowi02}
Richard~A. Johnson and Dean~W. Wichern.
\newblock \emph{Applied Multivariate Statistical Analysis}.
\newblock Prentice Hall, 5th edition, 2002.

\bibitem[Johnstone and Silverman(2004)]{josi02a}
I.~M. Johnstone and B.~W. Silverman.
\newblock Needles and straw in haystacks: {E}mpirical {B}ayes estimates of
  possibly sparse sequences.
\newblock \emph{Annals of Statistics}, 32:\penalty0 1594--1649, 2004.

\bibitem[Johnstone(2001)]{john00c}
Iain~M. Johnstone.
\newblock On the distribution of the largest eigenvalue in principal components
  analysis.
\newblock \emph{Annals of Statistics}, 29:\penalty0 295--327, 2001.

\bibitem[Johnstone(2006)]{john06}
Iain~M. Johnstone.
\newblock Canonical correlation analysis and {J}acobi ensembles: Tracy-{W}idom
  limits and rates of convergence.
\newblock Manuscript, 50pp, August 2006.

\bibitem[Johnstone and Lu(2004)]{jolu04}
Iain~M. Johnstone and Arthur~Yu Lu.
\newblock Sparse principal components analysis.
\newblock Technical report, Stanford University, Dept. of Statistics, 2004.
\newblock Tentatively accepted, {\it J.A.S.A.}

\bibitem[Jolliffe(2002)]{joll02}
I.~T. Jolliffe.
\newblock \emph{Principal Component Analysis}.
\newblock Springer, 2nd edition, 2002.

\bibitem[Koev(2006)]{koev06}
Plamen Koev.
\newblock Software {\tt mhg, mhgi} for hypergeometric function of a matrix
  argument.
\newblock \url{http://www-math.mit.edu/~plamen/}, 2006.

\bibitem[Koev and Edelman(2006)]{koed06}
Plamen Koev and Alan Edelman.
\newblock The efficient evaluation of the hypergeometric function of a matrix
  argument.
\newblock \emph{Math. Comp.}, 75\penalty0 (254):\penalty0 833--846
  (electronic), 2006.

\bibitem[Ledoit and Wolf(2004)]{ledo04}
Olivier Ledoit and Michael Wolf.
\newblock A well-conditioned estimator for large-dimensional covariance
  matrices.
\newblock \emph{Journal of Multivariate Analysis}, 88:\penalty0 365--411, 2004.

\bibitem[Lucas et~al.(2006)Lucas, Carvalho, Wang, Bild, Nevins, and
  West]{lcwm06}
Joe Lucas, Carlos Carvalho, Quanli Wang, Andrea Bild, Joe Nevins, and Mike
  West.
\newblock Sparse statistical modelling in gene expression genomics.
\newblock In K.~A. Do, P.~Mueller, and M.~Vannucci, editors, \emph{Bayesian
  Inference for Gene Expression and Proteomics}, pages 155--176. Cambridge
  University Press, 2006.

\bibitem[Mardia et~al.(1979)Mardia, Kent, and Bibby]{mkb79}
K.~V. Mardia, J.~T. Kent, and J.~M. Bibby.
\newblock \emph{Multivariate Analysis}.
\newblock Academic Press, 1979.

\bibitem[Mar\v{c}enko and Pastur(1967)]{mapa67}
V.~A. Mar\v{c}enko and L.~A. Pastur.
\newblock Distributions of eigenvalues of some sets of random matrices.
\newblock \emph{Math. USSR-Sb.}, 1:\penalty0 507--536, 1967.

\bibitem[Menozzi et~al.(1978)Menozzi, Piazza, and Cavalli-Sforza]{mpcs78}
P~Menozzi, A~Piazza, and L~Cavalli-Sforza.
\newblock {Synthetic maps of human gene frequencies in Europeans}.
\newblock \emph{Science}, 201\penalty0 (4358):\penalty0 786--792, 1978.

\bibitem[Mood(1951)]{mood51}
A.~M. Mood.
\newblock On the distribution of the characteristic roots of normal
  second-moment matrices.
\newblock \emph{The Annals of Mathematical Statistics}, 22:\penalty0 266--273,
  1951.

\bibitem[Muirhead(1982)]{muir82}
R.~J. Muirhead.
\newblock \emph{Aspects of Multivariate Statistical Theory}.
\newblock Wiley, 1982.

\bibitem[Onatski(2006)]{onat06b}
Alexei Onatski.
\newblock Asymptotic distribution of the principal components estimator of
  large factor models when factors are relatively weak.
\newblock URL \url{http://www.columbia.edu/~ao2027/papers1.html}.
\newblock Dept. of Economics, Columbia University, 2006.

\bibitem[Paul(2004)]{paul04}
Debashis Paul.
\newblock Asymptotics of sample eigenstructure for a large dimensional spiked
  covariance model.
\newblock Technical report, Department of Statistics, Stanford University,
  2004.
\newblock {\it Statistica Sinica}, to appear.

\bibitem[Pearson(1901)]{pear01}
K.~Pearson.
\newblock On lines and planes of closest fit to systems of points in space.
\newblock \emph{Philosophical Magazine}, 2\penalty0 (6):\penalty0 559--572,
  1901.

\bibitem[Persson(2002)]{pers02}
Per-Olof Persson.
\newblock Numerical methods for random matrices.
\newblock \url{http://www.mit.edu/~persson/numrand_report.pdf}, 2002.

\bibitem[Potters et~al.(2005)Potters, Bouchaud, and Laloux]{pbl05}
M.~Potters, J.~P. Bouchaud, and L.~Laloux.
\newblock Financial applications of random matrix theory: Old laces and new
  pieces, 2005.
\newblock URL
  \url{http://www.citebase.org/abstract?id=oai:arXiv.org:physics/0507111}.

\bibitem[Ram\'{i}rez et~al.(2006)Ram\'{i}rez, Rider, and Vir\'{a}g]{rrv06}
Jos\'{e} Ram\'{i}rez, Brian Rider, and B\'{a}lint Vir\'{a}g.
\newblock {Beta ensembles, stochastic Airy spectrum, and a diffusion}, 2006.

\bibitem[Reimann et~al.(1996)Reimann, Van~den Broeck, and Bex]{rebb96}
P.~Reimann, C.~Van~den Broeck, and G.~J. Bex.
\newblock A {G}aussian scenario for unsupervised learning.
\newblock \emph{Journal of Physics A: Mathematical and General}, 29\penalty0
  (13):\penalty0 3521--3535, 1996.

\bibitem[Roweis and Saul(2000)]{rosa00}
Sam~T. Roweis and Lawrence~K. Saul.
\newblock {Nonlinear Dimensionality Reduction by Locally Linear Embedding}.
\newblock \emph{Science}, 290\penalty0 (5500):\penalty0 2323--2326, 2000.

\bibitem[Roy(1939)]{roy39}
S.~N. Roy.
\newblock $p-$statistics or some generalizations in analysis of variance
  appropriate to multivariate problems.
\newblock \emph{Sankhy\={a}}, 4:\penalty0 381--396, 1939.

\bibitem[Soshnikov(2002)]{sosh01a}
Alexander Soshnikov.
\newblock A note on universality of the distribution of the largest eigenvalues
  in certain classes of sample covariance matrices.
\newblock \emph{J. Statistical Physics}, 108:\penalty0 1033--1056, 2002.

\bibitem[Soshnikov(2006)]{sosh06}
Alexander Soshnikov.
\newblock Poisson statistics for the largest eigenvalues in random matrix
  ensembles.
\newblock In \emph{Mathematical physics of quantum mechanics}, volume 690 of
  \emph{Lecture Notes in Phys.}, pages 351--364. Springer, Berlin, 2006.

\bibitem[Soshnikov and Fyodorov(2005)]{sofy05}
Alexander Soshnikov and Yan~V. Fyodorov.
\newblock On the largest singular values of random matrices with independent
  {C}auchy entries.
\newblock \emph{J. Math. Phys.}, 46\penalty0 (3):\penalty0 033302, 15, 2005.

\bibitem[Stein(ca. 1977)]{stei77}
C.~Stein.
\newblock Estimation of a covariance matrix.
\newblock Unpublished manuscript, Stanford University, ca. 1977.

\bibitem[Szeg\"{o}(1967)]{szeg67}
Gabor Szeg\"{o}.
\newblock \emph{Orthogonal Polynomials, 3rd edition}.
\newblock American Mathematical Society, 1967.

\bibitem[Tenenbaum et~al.(2000)Tenenbaum, Silva, and Langford]{tsl00}
Joshua~B. Tenenbaum, Vin~de Silva, and John~C. Langford.
\newblock {A Global Geometric Framework for Nonlinear Dimensionality
  Reduction}.
\newblock \emph{Science}, 290\penalty0 (5500):\penalty0 2319--2323, 2000.

\bibitem[Tracy and Widom(1994)]{trwi94}
Craig~A. Tracy and Harold Widom.
\newblock Level-spacing distributions and the {A}iry kernel.
\newblock \emph{Communications in Mathematical Physics}, 159:\penalty0
  151--174, 1994.

\bibitem[Tracy and Widom(1996)]{trwi96}
Craig~A. Tracy and Harold Widom.
\newblock On orthogonal and symplectic matrix ensembles.
\newblock \emph{Communications in Mathematical Physics}, 177:\penalty0
  727--754, 1996.

\bibitem[Tracy and Widom(1998)]{trwi98}
Craig~A. Tracy and Harold Widom.
\newblock Correlation functions, cluster functions, and spacing distributions
  for random matrices.
\newblock \emph{J. Statistical Physics}, 92:\penalty0 809--835, 1998.

\bibitem[Tulino and Verdu(2004)]{tuve04}
Antonia Tulino and Sergio Verdu.
\newblock \emph{Random Matrix Theory and Wireless Communications}.
\newblock Now Publishers Inc, 2004.

\bibitem[Wachter(1978)]{wach78}
Kenneth~W. Wachter.
\newblock The strong limits of random matrix spectra for sample matrices of
  independent elements.
\newblock \emph{Annals of Probability}, 6:\penalty0 1--18, 1978.

\bibitem[Wachter(1980)]{wach80}
Kenneth~W. Wachter.
\newblock The limiting empirical measure of multiple discriminant ratios.
\newblock \emph{Annals of Statistics}, 8:\penalty0 937--957, 1980.

\bibitem[Watkin and Nadal(1994)]{wana94}
T.~L.~H. Watkin and J.-P. Nadal.
\newblock Optimal unsupervised learning.
\newblock \emph{Journal of Physics A: Mathematical and General}, 27\penalty0
  (6):\penalty0 1899--1915, 1994.

\bibitem[Wigner(1955)]{wign55}
Eugene~P. Wigner.
\newblock Characteristic vectors of bordered matrices of infinite dimensions.
\newblock \emph{Annals of Mathematics}, 62:\penalty0 548--564, 1955.

\bibitem[Wigner(1958)]{wign58}
Eugene~P. Wigner.
\newblock On the distribution of the roots of certain symmetric matrices.
\newblock \emph{Annals of Mathematics}, 67:\penalty0 325--328, 1958.

\bibitem[Wishart(1928)]{wish28}
John Wishart.
\newblock The generalised product moment distribution in samples from a normal
  multivariate population.
\newblock \emph{Biometrika}, 20A\penalty0 (1/2):\penalty0 32--52, 1928.

\bibitem[Yang and Berger(1994)]{yabe94}
Ruoyang Yang and James~O. Berger.
\newblock Estimation of a covariance matrix using the reference prior.
\newblock \emph{Annals of Statistics}, 22:\penalty0 1195--1211, 1994.

\end{thebibliography}




 \frenchspacing








\end{document}